\documentclass[]{article}

\usepackage[utf8]{inputenc}
\usepackage[T1]{fontenc}

\usepackage[english]{babel}
\usepackage[margin=1in]{geometry}
\usepackage{etoolbox}
\usepackage[]{footmisc}

\usepackage[numbers]{natbib}
\usepackage{amsmath}
\usepackage{amsthm}
\usepackage{amssymb}

\usepackage{newtxtext}
\usepackage{newtxmath}

\usepackage{bbm}
\usepackage[inline]{enumitem}
\usepackage{array}

\usepackage{chngcntr}
\usepackage{mathrsfs}
\usepackage{mathtools}
\usepackage{graphicx}
\usepackage{scalerel}
\usepackage[nodisplayskipstretch]{setspace}
\usepackage[%
    colorlinks=true,
    breaklinks=true,
    bookmarks=true,
    urlcolor=blue,
    citecolor=blue,
    linkcolor=blue,
    bookmarksopen=false,
    draft=false%
]{hyperref}
\usepackage{microtype}

\setenumerate{label=(\roman*),itemsep=3pt,topsep=3pt}
\newcolumntype{P}[1]{>{\centering\arraybackslash}p{#1}}

\newcommand{\email}[1]{\href{mailto:#1}{\texttt{#1}}} 

\newcommand{\keywords}{\noindent \textbf{Keywords:}\ }

\renewcommand{\bar}{\overline}

\newcommand{\condition}[2]{#1~$\Longrightarrow$~#2}

\newcommand{\EE}{\mathbb{E}}

\newcommand{\RR}{\mathbb{R}}

\newcommand{\QQ}{\mathbb{Q}}

\newcommand{\NN}{\mathbb{N}}

\newcommand{\DD}{\mathcal{D}}

\newcommand{\KK}{\mathbb{K}}
\newcommand{\one}{\mathbbm{1}}
\newcommand{\inv}[1][1]{{}^{-#1}}
\newcommand{\cleq}{\preccurlyeq}

\renewcommand{\subseteq}{\subset}

\DeclareMathOperator{\Gr}{\mathrm{Gr}}
\DeclareMathOperator{\Dom}{\mathrm{Dom}}
\DeclareMathOperator*{\cl}{\mathrm{cl}}
\DeclareMathOperator*{\epi}{\mathrm{epi}}

\newcommand{\lsc}{\textit{lsc}}
\newcommand{\usc}{\textit{usc}}
\newcommand{\FPTlsc}{\textit{FPTlsc}}
\newcommand{\FPTusc}{\textit{FPTusc}}
\newcommand{\lminsc}{\textit{lmsc}}
\newcommand{\umaxsc}{\textit{umsc}}
\newcommand{\linfsc}{\textit{lisc}}
\newcommand{\usupsc}{\textit{ussc}}
\newcommand{\FPTlminsc}{\textit{FPTlmsc}}
\newcommand{\FPTlinfsc}{\textit{FPTlisc}}
\DeclarePairedDelimiter{\set}{\{}{\}}
\DeclarePairedDelimiter{\abs}{|}{|}

\DeclarePairedDelimiter{\paren}{(}{)}

\newcommand{\seq}{\paren}

\newtheorem{theorem}{Theorem}
\newtheorem{lemma}{Lemma}
\newtheorem{corollary}{Corollary}

\theoremstyle{definition}
\newtheorem{definition}{Definition}
\newtheorem{example}{Example}

\theoremstyle{remark}
\newtheorem{remark}{Remark}

\makeatletter
\let\@fnsymbol\@arabic
\makeatother

\title{%
	Continuity of Parametric Optima for Possibly Discontinuous Functions and Noncompact
	Decision Sets
}

\date{%
	September 10, 2021
}

\author{
	Eugene A. Feinberg\thanks{Department of Applied Mathematics and Statistics, State University of New York at Stony Brook, Stony Brook, NY 11794-3600, USA, %
		\email{eugene.feinberg@stonybrook.edu}} \and
	Pavlo O. Kasyanov\thanks{Institute for Applied System Analysis, National Technical University of Ukraine ``Igor Sikorsky Kyiv Polytechnic Institute", Peremogy ave., 37, build, 35, 03056, Kyiv, Ukraine, %
		\email{kasyanov@i.ua}} \and
	David N. Kraemer\thanks{Department of Applied Mathematics and Statistics, State University of New York at Stony Brook, Stony Brook, NY 11794-3600, USA, %
		\email{david.kraemer@stonybrook.edu}}
}

\begin{document}

\maketitle

\onehalfspacing

\begin{abstract}
	This paper investigates continuity properties of value functions and solutions for parametric optimization problems. These problems are important in operations research, control, and economics because optimality equations are their particular cases. The classic fact, Berge's maximum theorem, gives sufficient conditions for continuity of value functions and upper semicontinuity of solution multifunctions.  Berge's maximum theorem assumes that the objective function is continuous and the multifunction of feasible sets is compact-valued. These assumptions are not satisfied in many applied problems, which historically has limited the relevance of the theorem. This paper generalizes Berge’s maximum theorem in three directions: (i) the objective function may not be continuous, (ii) the multifunction of feasible sets may not be compact-valued, and (iii) necessary and sufficient conditions are provided. To illustrate the main theorem, this paper provides applications to inventory control and to the analysis of robust optimization over possibly noncompact action sets and discontinuous objective functions.
	
	\keywords{%
	    Berge's maximum theorem;
	    continuity;
	    multifunction;
	    inventory control;
	    robust optimization
	}
\end{abstract}

\section{Introduction}
\label{sec:introduction}

This paper studies continuity properties of value functions and solution multifunctions for parametric optimization problems with possibly discontinuous objective functions and noncompact decision sets.  The classic result that states continuity of the value function under general conditions is Berge's maximum theorem~\cite[Chapter 6, Section 3]{berge1963topological}. We consider the minimization problem $(X, Y, \Phi, u)$
defined by the equation
\begin{equation}
    u^*(x) = \inf_{y \in \Phi(x)} u(x,y),\qquad\qquad x\in X,
    \label{eq:value-fn}
\end{equation}
for topological spaces $X$ and $Y$, a multifunction $\Phi : X \to 2^Y$, and objective function $u : X \times Y \to \bar{\RR}$, where $\inf \emptyset = +\infty$ by convention, $\bar{\RR} = \RR \cup \set{\pm \infty}$, and $\RR$ is the real line.  Berge's maximum theorem states that, if $u$ is continuous and $\Phi$ is continuous with compact values, then the value function $u^*$ is continuous, and the solutions multifunction
\begin{equation}
    \Phi^*(x) = \set{y \in \Phi(x) : u(x,y) = u^*(x)}
    \label{eq:solutions-multifn}
\end{equation}
is upper semicontinuous (\usc{}) with nonempty and compact values.

Berge's maximum theorem is a corollary of three statements guaranteeing the following properties of the value function and solution multifunction: upper semicontinuity of $u^*,$ lower semicontinuity of $u^*,$ and upper semicontinuity of $\Phi^*$ and non-emptiness and compactness of the sets $\Phi^*(x);$ see ~\citet[pp.~569--570]{aliprantis1999infinite},~\citet[pp.~115--116]{berge1963topological}, or~\citet[pp.~82--84]{hu1997handbook}.  In particular, the function $u^*$ is \usc{}, if the function $u$ is \usc{}, and if the multifunction $\Phi$ is \lsc{}.  The function $u^*$ is lower semicontinuous (\lsc{}), if the function $u$ is \lsc{}, and if the multifunction $\Phi$ is \usc{} and compact-valued.

Continuity of value functions is an important consideration for many applications. In reinforcement learning, many algorithms follow the paradigm of value function approximation, including \textsc{sarsa}, deep Q-learning, and actor-critic methods; see~\citet{bertsekas2019reinforcement},~\citet{sutton2018reinforcement}, and~\citet{szepesvari2010algorithms} for general overviews of this paradigm. Neural function approximation architecture has become ubiquitous because such models can efficiently estimate continuous functions. This is a consequence of several fundamental theorems of analysis, including the Stone-Weierstrass theorem (see, for example,~\citet{cybenko1989approximation} for the classic result) and the Kolmogorov-Arnold superposition theorem (see, for example~\citet{montanelli2020error} as a recent development). In economics, continuity of the value function (also known as the marginal function) is a principal objective in the analysis of consumer behavior and economic planning; see, for example,~\citet[\S 7.3]{varian1992microeconomic} in the context of microeconomic agents and~\citet[Theorems 3.6 and Theorem 4.6]{stokey1989recursive} for macroeconomic theory.

This paper generalizes Berge's maximum theorem by relaxing the assumptions that the objective function $u$ is continuous, real-valued, and the multifunction $\Phi$ is upper semicontinuous and compact-valued. We study general conditions that ensure upper semicontinuity and lower semicontinuity of $u^*$ as well as the upper semicontinuity of $\Phi^*$, including a number of necessary and sufficient conditions. There are several results obtained in this direction after Berge's maximum theorem was published.  \citet{tian1992maximum} introduced the property of feasible path transfer upper semicontinuity, showed that it is necessary and sufficient for upper semicontinuity of $u^*,$ and generalized Berge's maximum theorem by using this property.  \citet{feinberg2013berge} introduced the notion of $\KK$-inf-compact functions and established sufficient conditions for $u^*$ to be \lsc{} when the sets $\Phi(x)$ are possibly noncompact. The results in \cite{feinberg2013berge} cover mostly metric spaces $X$ and $Y.$ \citet{feinberg2014berge} introduced the notion of $\KK\NN$-inf-compact functions for Hausdorff metric spaces and extended Berge's maximum theorem to noncompact sets $\Phi(x).$ \citet[Proposition 4.4]{bonnans2000perturbation} introduced local Berge's maximum theorem when $X$ is a Banach space, and the sets $\Phi(x)$ are possibly noncompact. \citet{feinberg2015continuity} described a general version of a local Berge's maximum theorem for Hausdorff topological spaces $X$ and $Y.$ \citet[Theorems 1.17, 7.41]{rockafellar2009variational} developed, as an application of their study of epi-continuous functions, a version of Berge's maximum theorem for the case when $X$ and $Y$ are Euclidean spaces and $\Phi(x) = Y$ for all $x \in X$.

This paper describes generalizations of Berge's maximum theorem and relevant results by combining the properties of feasible path transfer upper semicontinuity, a new condition called lower min-semicontinuity, and $\KK\NN$-inf-compactness of objective functions. Section~\ref{sec:background} presents a summary of relevant background and notation. Section~\ref{sec:minimization} gives necessary and sufficient conditions for upper and lower semicontinuity of the value function as well as for existence and upper semicontinuity of solutions.  It formulates natural generalizations of the results in \citet{tian1992maximum}, \citet{feinberg2013berge, feinberg2014berge}, and \citet{feinberg2015continuity}.  These are also natural generalizations of \citet[Theorems 1.17, 7.41]{rockafellar2009variational} to topological spaces $X$ and $Y$ and to multifunctions $\Phi(x)$ depending on $x\in X$.  The general form of Berge's maximum theorem is Theorem~\ref{thm:general-local-max} stating that feasible path transfer upper semicontinuity and $\KK\NN$-inf-compactness imply the conclusions of Berge's maximum theorem.  Lower min-semicontinuity is introduced in Definition~\ref{def:lminsc} as a necessary and sufficient condition for $u^*$ to be \lsc{} and $\Phi^*$ to be nonempty, as stated in Theorem~\ref{thm:lminsc}. A weaker property, lower inf-semicontinuity, is introduced in Definition~\ref{def:linfsc} and is shown to be equivalent to \lsc{} of $u^*$. The local versions of Berge's maximum theorem, presented in Theorems~\ref{thm:main-berge-max} and \ref{thm:general-local-max}, are more general results than the corresponding theorems in both~\citet[Theorem 7]{feinberg2015continuity} and~\citet[Theorem 1]{tian1992maximum}. Section~\ref{subsection:formulations} compares three different formulations of the parametric optimization problem: the case where $\Phi(x) \equiv Y$, the case where $u < +\infty$ on $\Gr(\Phi)$, and the formulation adopted by Equation~(\ref{eq:value-fn}). The interpretation of the formulation with $\Phi(x) \equiv Y$ is parametric optimization without explicit constraints, and this formulation is studied in~\citet{rockafellar2009variational} for the case $Y = \RR^n$.  The formulations are shown to be mostly equivalent with respect to determining feasible path transfer upper semicontinuity, lower min-semicontinuity, and lower inf-semicontinuity. Section~\ref{sec:minimization-metric} discusses the general parametric minimization results of Section~\ref{sec:minimization} in the context of metric spaces and provides some stronger results for this setting. Section~\ref{sec:examples} provides a number of examples employing and differentiating between the properties studied in Section~\ref{sec:minimization}. It also includes applications to inventory control with possibly bounded orders, possibly finite storage capacity, and with or without backorders. Section~\ref{sec:robust-optimization} considers applications to minimax equations for robust optimization, including results on continuity of minimax functions, existence of equilibria, and continuity of equilibria.

In conclusion of this section, we would like to mention that there are several extensions of Berge's maximum theorem to various additional settings. \citet{debreu1967neighboring} states the results for parametric optimization of non-numerical functions. Generalizations for the setting in which $X$ and $Y$ are vector spaces and $u$ is convex include~\citet{bank1982nonlinear},~\citet{sundaram1996first}, and~\citet{bertsekas2007set}, and, when $u$ is quasiconvex, \citet{terazono2015continuity}. Additional references can also be found in~\citet[Section 1]{feinberg2014berge}. 

\section{Notation} \label{sec:background}

Let $\Phi : X \to 2^Y$ be a multifunction; that is, $\Phi(x) \subseteq Y$ for each $x \in X$.  On a subset $E \subseteq X$, we write the graph of $\Phi$ on $E$ as
\begin{equation*}
  \Gr_E(\Phi) = \set{(x,y) \in E \times Y : y \in \Phi(x)},
\end{equation*}
where $E \subseteq X$ (we omit the subscript when $E = X$), and domain
\begin{align*}
  \Dom(\Phi) = \set{x \in X : \Phi(x) \ne \emptyset}.
\end{align*}
When $\Dom(\Phi) = X$, we call $\Phi$ a \emph{strict} multifunction.  Let $x \in \Dom(\Phi)$ be given. A \emph{neighborhood} of $x$ is an open set $U \subset X$ containing $x$, and we write $U(x)$ to express the dependency on $x$ explicitly. The multifunction $\Phi$ is \emph{upper semicontinuous} (\usc{}) at $x$ if, for each open set $V \subseteq Y$ with $\Phi(x) \subset V$, there exists a neighborhood $U(x)$ such that $\Phi(x') \subseteq V$ for each $x' \in U(x)$.  The multifunction $\Phi$ is \emph{lower semicontinuous} (\lsc{}) at $x$ if for each open set $V \subseteq Y$ with $\Phi(x) \cap V \ne \emptyset$, there exists a neighborhood $U(x)$ such that $\Phi(x') \cap V \ne \emptyset$ for each $x' \in U(x)$. If $\Phi$ is both \usc{} and \lsc{} at $x$, it is said to be \emph{continuous at $x$}. If $\Phi$ is \usc{} (resp., \lsc{}, continuous) at $x$ for each $x \in \Dom(\Phi)$, then $\Phi$ is \usc{} (resp., \lsc{}, continuous). For a set $E \subset \Dom(\Phi)$, we write $f : E \subseteq X \to Y$ to indicate a function with values in the space $Y$ defined on a subset $E$ of the space $X$.  When $E = X$, we keep the usual notation $f : X \to Y$. A \emph{selector} for $\Phi$ on $E$ is a function $\sigma : E \subseteq X \to Y$ such that $\sigma(x) \in \Phi(x)$ for each $x \in E$. We denote the set of selectors by $\Sigma(\Phi;E)$ and use the shorthand $\Sigma(\Phi) := \Sigma(\Phi; \Dom(\Phi))$.

We always consider the order topology on the set of real numbers $\RR$ and on the extended set of real numbers $\bar{\RR} = \RR \cup \set{\pm\infty}.$ A function $f : E \subseteq X \to \bar{\RR}$ is \emph{lower semicontinuous} (\emph{\lsc{}}) \emph{at $x \in E$} if, each net $\seq{x_\alpha}_{\alpha \in I}$ with values in $E$ and with $x_\alpha \to x$ satisfies $f(x) \leq \liminf_{\alpha \in I} f(x_\alpha)$, and $f$ is \emph{upper semicontinuous} (\emph{\usc{}}) \emph{at $x \in E$} if $-f$ is \lsc{} at $x$; see, e.g.,~\citet[p. 43]{aliprantis1999infinite} for an approach to this definition.  Of course, $f$ is continuous at $x \in E$ if and only if $f$ is \lsc{} and \usc{} at $x$. In this case, for every net $\seq{x_\alpha}_{\alpha \in I}$ such that $x_\alpha \to x$, we have $\lim_{\alpha \in I} f(x_\alpha) =\liminf_{\alpha \in I} f(x_\alpha) = \limsup_{\alpha \in I} f(x_\alpha) = f(x)$, and we write $f(x_\alpha) \to f(x)$ as $x_\alpha \to x$.  If $f$ is \lsc{} (resp., \usc{}, continuous) at $x$ for each $x \in E$, then $f$ is \emph{\lsc{}} (resp., \emph{\usc{}}, \emph{continuous}) \emph{on $E$}. The function $f$ is \lsc{} at $x$ if and only if for each $\lambda \in \RR$, if $\lambda < f(x)$, there exists a neighborhood $U(x)$ such that each $x' \in U(x)$ satisfies $\lambda < f(x')$; see, e.g., Definition 1 and Proposition 3 in~\citet[Chapter 4, Section 6.2]{bourbaki2013general}, where Proposition 3 deals with filters, an equivalent approach to nets. The epigraphical mapping of the function $f$, denoted $\epi f$, is defined as $(\epi f)(x) = \set{\lambda \in \RR : f(x) \leq \lambda}$. For $\lambda \in \RR$, we denote the level set of $f$ by $\mathcal{D}_f(\lambda ; E) = \set{x \in E : f(x) \leq \lambda}$. The indicator function of the set $E$ is denoted as $\one_E: X \to \set{0,1}$, meaning $\one_E(x) = 1$ iff $x \in E$.

\section{Continuity of minima} \label{sec:minimization}

For topological spaces $X$ and $Y$, a multifunction $\Phi : X \to 2^Y$, and objective function $u : \Gr(\Phi) \subseteq X \times Y \to \bar\RR$, we consider the minimization problem defined by the tuple $(X, Y, \Phi, u)$  with the \emph{value function} $u^* : \Dom(\Phi) \subseteq X \to \bar\RR$ and with the \emph{solutions multifunction} $\Phi^* : \Dom(\Phi) \subset X \to 2^Y$ defined in~(\ref{eq:value-fn}) and~(\ref{eq:solutions-multifn}). This section investigates conditions on $\Phi$ and $u$ which guarantee basic continuity properties of the value function $u^*$ and the solutions multifunction $\Phi^*$.

The section is divided into four subsections. Section~\ref{subsection:usc-value} considers necessary and sufficient conditions for the value function $u^*$ to be \usc{}; Section~\ref{subsection:lsc-value} considers conditions for $u^*$ to be \lsc{}; Section~\ref{subsection:main-results} considers conditions for the solutions $\Phi^*$ to be nonempty, \usc{}, and compact-valued, and it states Theorem~\ref{thm:main-berge-max}, the main result of this section, which generalizes Theorem 7 in ~\citet{feinberg2015continuity} and Theorem 1 in~\citet{tian1992maximum}. Section~\ref{subsection:formulations} considers the equivalence between different formulations of the parametric optimization problem defined in equations (\ref{eq:value-fn}) and (\ref{eq:solutions-multifn}). Theorems~\ref{thm:main-berge-max} and \ref{thm:general-local-max} present more general formulations of the local Berge's maximum theorem than any known formulation.

\subsection{Upper semicontinuity of values}
\label{subsection:usc-value}

\citet{tian1992maximum} introduced the notion of feasible path transfer semicontinuity for real-valued objective functions $u$. The following definition slightly extends this notion to include cases when $u$ can take infinite values.

\begin{definition}
    [{Feasible path transfer semicontinuity; cp.~\citet[Definition 1]{tian1992maximum}}]
    Let $E \subset \Dom(\Phi)$.  The function $u : \Gr(\Phi) \subseteq X \times Y \to \bar\RR$ is called \emph{feasible path transfer upper semicontinuous} (\emph{\FPTusc{}}) on $\Gr_{E}(\Phi)$ if for each $x \in E$, $y \in\Phi(x)$, and $\gamma \in \RR$ with $u(x,y) < \gamma$, there exists a neighborhood $U(x)$ such that for each $x' \in U(x)$ there exists $y' \in \Phi(x')$ such that $u(x',y') < \gamma$.  The function $u$ is \emph{feasible path transfer lower semicontinuous} (\emph{\FPTlsc{}}) on $\Gr_E(\Phi)$ if $-u$ is \FPTusc{} on $\Gr_E(\Phi)$.
    \label{def:fptusc}
\end{definition}

Alternatively, one can reformulate Definition~\ref{def:fptusc} using selectors. The function $u : \Gr(\Phi) \subseteq X \times Y \to \bar\RR$ is \FPTusc{} on $\Gr_E(\Phi)$ if and only if for each $x \in E$, $y \in \Phi(x)$, and $\gamma \in \RR$ with $u(x,y) < \gamma$, there is a neighborhood $U(x)$ and a selector $y' \in \Sigma(\Phi ; U(x))$ such that $u(x', y'(x')) < \gamma$ for each $x' \in U(x)$.  Such a selector $y'$ is called a \emph{feasible path}.  According to Proposition 1 in~\citet{tian1992maximum}, the value function $u^*$ is \usc{} if and only if the objective function $u$ is \FPTusc{} on $\Gr_E(\Phi)$.

The following theorem, due to~\citet{tian1992maximum}, establishes necessary and sufficient conditions for the function $u^*$ to be \usc{}. In particular, feasible path transfer upper semicontinuity of $u$ completely characterizes upper semicontinuity of $u^*$. There are minor technical differences between Theorem~\ref{thm:general-usc-value} and \cite[Proposition 1]{tian1992maximum}: Theorem~\ref{thm:general-usc-value} is a local statement, and it uses Definition~\ref{def:fptusc} allowing infinite values of $u.$ However, the proof of \cite[Proposition 1]{tian1992maximum} applies to Theorem~\ref{thm:general-usc-value}.

\begin{theorem}
    [{cf.~\citet[Proposition 1]{tian1992maximum}}]
    The function $u^*: \Dom(\Phi) \subset X \to \bar\RR$ is \usc{} at $x \in \Dom(\Phi)$ if and only if the function $u : \Gr(\Phi) \subset X \times Y \to \bar \RR$ is \FPTusc{} on $\Gr_{\set{x}}(\Phi)$.
    \label{thm:general-usc-value}
\end{theorem}

The global version of Theorem~\ref{thm:general-usc-value} is given in the next corollary.

\begin{corollary}
    [{cf.~\citet[Proposition 1]{tian1992maximum}}]
    The function $u^*: \Dom(\Phi) \subset X \to \bar\RR$ is \usc{} if and only if the function $u : \Gr(\Phi) \subset X \times Y \to \bar \RR$ is \FPTusc{} on $\Gr(\Phi)$.
    \label{cor:general-usc-value}
\end{corollary}

The next theorem provides a number of sufficient conditions for the function $u^*$ to be \usc{}. Condition~\ref{item:berge-usc-value} corresponds to the classic sufficient condition for the value function $u^*$ to be \usc{}; see~\citet[Chapter 6, Section 3, Theorem 1]{berge1963topological} or~\citet[Lemma 17.29]{aliprantis1999infinite}. The proof for condition~\ref{item:berge-usc-value} in Theorem~\ref{thm:sufficient-usc} is given because it includes extended real-valued objective functions $u$. Condition~\ref{item:rockafellar-usc-value} in Theorem~\ref{thm:sufficient-usc} is a generalization of the similar result for Euclidean spaces in~\citet[Theorem 1.17(c)]{rockafellar2009variational} to general topological spaces.  We note that $u^*$ is automatically \usc{} on each $x \in X \setminus \Dom(\Phi)$, because $u^*$ remains unchanged when $u$ is extended from $\Gr(\Phi)$ to $X \times Y$ by setting $u(z,y) := -\infty$ if $y \in Y \setminus \Phi(z)$ for each $z \in X$. This is a consequence of the reduction to the case described in condition~\ref{item:rockafellar-usc-value} of Theorem~\ref{thm:sufficient-usc}. Condition~\ref{item:fptusc-nets} of Theorem~\ref{thm:sufficient-usc} is a novel sufficient condition inspired by epi-upper semicontinuity, a concept developed in~\citet[Definition 7.39, Theorem 7.41]{rockafellar2009variational} for Euclidean spaces. In metric spaces, replacing the nets in Condition~\ref{item:fptusc-nets} with sequences yields a necessary and sufficient condition for $u$ to be \FPTusc{}; see Theorem~\ref{thm:metric-fptusc}, below, for details.

\begin{theorem}
    The following conditions are sufficient for the function $u^*: \Dom(\Phi) \subset X \to \bar\RR$ to be \usc{} at $x \in \Dom(\Phi)$,
    \begin{enumerate}
        \item \label{item:berge-usc-value} the function $u : \Gr(\Phi) \subset X \times Y \to \bar\RR$ is \usc{} on $\Gr_{\set{x}}(\Phi)$ and the multifunction $\Phi$ is \lsc{} on $\Gr_{\set{x}}(\Phi)$;
        \item \label{item:rockafellar-usc-value} for each $y \in \Phi(x)$ the function $u(\:\cdot\:,y)$ is \usc{} at $x$ and there is a neighborhood $U(x)$ such that $\Phi(x') = A \subset Y$ for each $x' \in U(x)$;
        \item \label{item:fptusc-nets} for each $y \in \Phi(x)$ and for each net $\seq{x_\alpha}_{\alpha \in I} \subset \Dom(\Phi)$ with $x_{\alpha} \to x$, there exists a net $\seq{y_\alpha}_{\alpha \in I}$ with $y_\alpha \in \Phi(x_\alpha)$ for each $\alpha \in I$ with $\limsup_{\alpha \in I} u(x_\alpha, y_\alpha) \leq u(x,y)$.
    \end{enumerate}
    \label{thm:sufficient-usc}
\end{theorem}

\begin{proof}
To prove that $u^*$ is \usc{} at $x$, let $\gamma \in \RR$.
\begin{enumerate}[%
    align=left,
    itemindent=2\parindent,
    labelwidth=\parindent,
    listparindent=\parindent,
    leftmargin=0pt
  ]
    \item If $u^*(x) < \gamma$, fix $y \in \Phi(x)$ such that $u(x,y) < \gamma$. Since $u$ is \usc{} on $\Gr_{\set{x}}(\Phi)$, there are neighborhoods $U_1(x)$ and $V(y)$ with $u(x',y') < \gamma$ for each $x' \in U_1(x)$ and $y' \in V(y)$. Since $\Phi$ is \lsc{} at $x$, there is a neighborhood $U_2(x)$ such that $\Phi(x') \cap V(y) \ne \emptyset$ for each $x' \in U_2(x)$. Then for each $x' \in U_1(x) \cap U_2(x)$ there is $y' \in \Phi(x') \cap V(y)$, such that $u^*(x') \leq u(x',y') < \gamma$, implying that $u^*$ is \usc{}.
    \item If $u^*(x) < \gamma$, fix $y \in \Phi(x)$ such that $u(x,y) < \gamma$. Since $u(\:\cdot\:,y)$ is \usc{} at $x$, there is a neighborhood $U_1(x)$ such that $u(x',y) < \gamma$ for each $x' \in U_1(x)$. Furthermore, there is a neighborhood $U_2(x)$ such that $y \in \Phi(x')$ for each $x' \in U_2(x)$. Then for each $x' \in U_1(x) \cap U_2(x)$, we have $y \in \Phi(x')$ and $u^*(x') \leq u(x',y) < \gamma$, implying that $u^*$ is \usc{}.
    \item Let the net $\seq{x_\alpha}_{\alpha \in I} \subset \Dom(\Phi)$ with $x_\alpha \to x$, and let $\varepsilon > 0$. Fix $y \in \Phi(x)$ such that $u(x,y) \leq u^*(x) + \varepsilon$. Then there is a net $y_\alpha \in \Phi(x_\alpha)$ for each $\alpha \in I$ such that $\limsup_{\alpha \in I} u^*(x_\alpha) \leq \limsup_{\alpha \in I} u(x_\alpha, y_\alpha) \leq u(x,y) \leq u^*(x) + \varepsilon$. Since $\varepsilon > 0$ was arbitrary, $u^*$ is \usc{}. \qedhere
\end{enumerate}
\end{proof}

\subsection{Lower semicontinuity of values and existence, compactness of solutions}
\label{subsection:lsc-value}

In this section, we consider necessary and sufficient conditions for the value function $u^*$ to be \lsc{}, the existence of solutions, and general sufficient conditions for the compactness of the solution sets $\Phi^*$.  We introduce the related notions of lower min-semicontinuity and lower inf-semicontinuity. A lower inf-semicontinuous function $u$ is necessary and sufficient for lower semicontinuity of $u^*$, whereas a lower min-semicontinuous function $u$ is necessary and sufficient to lower semicontinuity of $u^*$ and existence of solutions $y^* \in \Phi^*(x)$. In addition, we extend to general topological spaces and extended-valued objective functions the results for $\KK\NN$-inf-compact functions, introduced in~\citet[Definition 1.3]{feinberg2014berge} and~\citet[Definition
1]{feinberg2015continuity} as a natural sufficient condition for $u^*$ to be \lsc{} and for $\Phi^*(x)$ to be nonempty and compact.

\begin{definition}
    [Lower min-semicontinuity]
    The function $u : \Gr(\Phi) \subset X \times Y \to \bar\RR$ is called \emph{lower min-semicontinuous} (\emph{\lminsc}) at $x \in \Dom(\Phi)$ if there exists $y^* \in \Phi(x)$ such that for every $\gamma \in \RR$ with $u(x,y^*) > \gamma$, there exists a neighborhood $U(x)$ such that $u(x',y') > \gamma$ for each $x' \in U(x)$ and $y' \in \Phi(x')$. The function $u$ is \emph{upper max-semicontinuous} (\emph{\umaxsc}) at $x \in \Dom(\Phi)$ if $-u$ is \lminsc{} at $x$.
    \label{def:lminsc}
\end{definition}

An example of a \lminsc{} function is $u : [0,1] \times [0,1] \to \RR$ given by the formula $u(x,y) = \one_{\QQ}(x-y)$, where for each $x \in [0,1]$ we choose $y^* \in [0,1]$ such that $x-y^* \notin \QQ$ and hence $u(x,y^*) = 0$. Note that this $u$ is nowhere \lsc{}. For a multifunction $\Gamma : \Gr(\Phi) \subset X \times Y \to 2^Z$, let $P_X^\Phi \Gamma$ denote the projection operator onto the space $X$ defined by $(P_X^\Phi \Gamma)(x) = \bigcup_{y \in \Phi(x)} \Gamma(x,y)$. We remark that the equality $\cl P_X^\Phi \epi u = \epi u^*$ holds in general for any function $u$. This is because the following statements are equivalent: (i) $\lambda \in (\cl P_X^\Phi \epi u)(x)$; (ii) there is a sequence $\lambda_n \in (P_X^\Phi \epi u)(x)$ with $\lambda_n \to \lambda$; (iii) there is a sequence $y_n \in \Phi(x)$ with $u(x,y_n) \leq \lambda_n$; (iv) $\lambda \in (\epi u^*)(x)$. As follows from the definitions, if $u$ is \lminsc{}, then $P_X^\Phi \epi u =\epi u^*$. This is the epigraphical projection equality as discussed in~\citet[Proposition 1.18]{rockafellar2009variational}. A weaker property given in the next definition is also useful for characterizing the value function $u^*$.

\begin{definition}
    [Lower inf-semicontinuity]
    The function $u : \Gr(\Phi) \subset X \times Y \to \bar\RR$ is called \emph{lower inf-semicontinuous} (\emph{\linfsc}) at $x \in \Dom(\Phi)$ if for all $\gamma \in \RR$ such that $u(x,y) > \gamma$ for each $y \in \Phi(x)$, there exists a neighborhood $U(x)$ such that $u(x',y') > \gamma$ for each $x' \in U(x)$ and $y' \in \Phi(x')$. The function $u$ is \emph{upper sup-semicontinuous} (\emph{\usupsc}) at $x \in \Dom(\Phi)$ if $-u$ is \linfsc{} at $x$.
    \label{def:linfsc}
\end{definition}

Every \lminsc{} function is also \linfsc{}, but the reverse implication is not true in general. A simple example of a \linfsc{} function that is not \lminsc{} is $u: [0,1] \times [0,1] \to \RR$ given by the formula $u(x,y) = y + \one_\QQ(x-y) + \one_{\set{0}}(y)$; see Example~\ref{ex:linfsc-not-lminsc} for details. If $u$ is \lminsc{} at $x$, then there exists a solution $y^* \in \Phi^*(x)$, but the solutions multifunction $\Phi^*(x)$ may be empty if $u$ is only \linfsc{}. The following theorem provides necessary and sufficient conditions for $u$ to be \linfsc{}. Condition~\ref{item:linfsc-lsc-min} of Theorem~\ref{thm:linfsc} states that lower inf-semicontinuity of $u$ is equivalent to lower semicontinuity of $u^*$.

\begin{theorem}
    [Necessary and sufficient conditions for {\linfsc{}}]
    For $x \in \Dom(\Phi)$, the following statements are equivalent:
    \begin{enumerate}
        \item The function $u : \Gr(\Phi) \subset X \times Y \to \bar \RR$ is \linfsc{} at $x$. \label{item:linfsc-def}
        \item For  $\lambda \in \RR$, and net $\seq{(x_\alpha,y_\alpha)}_{\alpha \in I} \subset \Gr(\Phi)$ with $x_\alpha \to x$, if $\sup_{\alpha \in I} u(x_\alpha, y_\alpha) < \lambda$, there exists $y \in \Phi(x)$ such that $u(x,y) < \lambda$. \label{item:linfsc-nets}
        \item The function $u^* : \Dom(\Phi) \subset X \to \bar\RR$ is \lsc{} at $x$. \label{item:linfsc-lsc-min}
    \end{enumerate}
    \label{thm:linfsc}
\end{theorem}

\begin{proof}
\condition{\ref{item:linfsc-def}}{\ref{item:linfsc-lsc-min}} Suppose $u^*(x) > \gamma$ for some $\gamma \in \RR$, and let $\gamma' \in (\gamma, u^*(x))$. Since $\gamma' < u(x,y)$ for each $y \in \Phi(x)$, there is a neighborhood $U(x)$ such that $u(x',y') > \gamma'$ for each $x' \in U(x)$ and $y' \in \Phi(x')$. Thus $u^*(x') \geq \gamma' > \gamma$ for each $x' \in U(x)$, so $u^*$ is \lsc{} at $x$.

\condition{\ref{item:linfsc-lsc-min}}{\ref{item:linfsc-nets}} Let $\lambda \in \RR$, and the net $\seq{(x_\alpha,y_\alpha)}_{\alpha \in I} \subset \Gr(\Phi)$ with $x_\alpha \to x$ be arbitrary. Suppose $\sup_{\alpha \in I} u(x_\alpha, y_\alpha) < \lambda$. Then $u^*(x) \leq \liminf_{\alpha \in I} u^*(x_\alpha) \leq \liminf_{\alpha \in I} u(x_\alpha, y_\alpha) \leq \sup_{\alpha \in I} u(x_\alpha, y_\alpha) < \lambda$, so there exists $y \in \Phi(x)$ such that $u(x,y) < \lambda$, as needed.

\condition{\ref{item:linfsc-nets}}{\ref{item:linfsc-def}} Suppose $u$ is not \linfsc{} at $x$. Then there is $y \in \Phi(x)$ and $\gamma \in \RR$ with $\gamma < u(x,y)$ such that, for every neighborhood $U(x)$, there are $x' \in U(x)$ and $y' \in \Phi(x')$ with $\gamma \geq u(x', y')$. Let $\gamma' \in (\gamma, u(x,y))$, and take the directed set $I = \set{U \subset X : \text{$U$ is a neighborhood of $x$}}$ with $U \cleq V$ iff $V \subset U$, it follows that there exists a net $\seq{(x_\alpha,y_\alpha)}_{\alpha \in I} \subset \Gr_{\set{x}}(\Phi)$ such that $\gamma \geq \sup_{\alpha \in I} u(x_\alpha, y_\alpha)$. Hence, $u(x,y) > \gamma' > \gamma \geq \sup_{\alpha \in I} u(x_\alpha, y_\alpha)$, so Condition~\ref{item:linfsc-nets} fails with $\gamma'$. The result follows by contrapositive.
\end{proof}

Theorem~\ref{thm:linfsc} is local in the sense that it characterizes lower semicontinuity at a single point $x \in X$. Since a function is \lsc{} on $E$ if and only if it is \lsc{} at each $x \in E$, the global version of Theorem~\ref{thm:linfsc} is an immediate corollary. In the global version, however, lower semicontinuity of $u^*$ is further characterized by closedness of $\Gr(\epi u^*)$.

\begin{corollary}
    Suppose $\Phi : X \to 2^Y$ is a strict multifunction. Then the following statements are equivalent:
    \begin{enumerate}
        \item The function $u : \Gr(\Phi) \subset X \times Y \to \bar \RR$ is \linfsc{} at each $x \in X$. \label{item:linfsc-def-2}
        \item The multifunction $\epi u^* : X \to 2^{\RR}$ has a closed graph. \label{item:linfsc-epi}
        \item The function $u^* : X \to \bar\RR$ is \lsc{} on $X$. \label{item:linfsc-lsc-min-2}
    \end{enumerate}
    \label{cor:linfsc}
\end{corollary}

\begin{proof}
Conditions~\ref{item:linfsc-def-2} and~\ref{item:linfsc-lsc-min-2} are equivalent in view of Theorem~\ref{thm:linfsc}. Conditions~\ref{item:linfsc-lsc-min-2} and~\ref{item:linfsc-epi} are equivalent, since $u^*$ is \lsc{} if and only if $\epi u^*$ has a closed graph.
\end{proof}

Although Corollary~\ref{cor:linfsc} assumes that $\Phi$ is strict, there are two ways to apply it to a possibly non-strict $\Phi$. First, Corollary~\ref{cor:linfsc} can be recovered by replacing $X$ with $\Dom(\Phi)$. In particular, for Condition~\ref{item:linfsc-epi}, this modification means that $\epi u^*$ has a closed graph in the subspace topology of $\Dom(\Phi) \times Y$. Second, the original minimization problem can be reduced to the problem with the same value function and with the multifunction identically equal to $Y$ and with setting $u(x,y)$ equal to $+\infty$ for $y\in Y\setminus \Phi(x);$ see Corollary~\ref{cor:strict-phi-bar-linfsc}.

Whereas lower inf-semicontinuity of $u$ is characterized by the property that $u^*$ is \lsc{}, the following theorem provides necessary and sufficient conditions for $u$ to be \lminsc{}. In particular, Condition~\ref{item:lminsc-lsc-min} states that $u$ is \lminsc{} at $x$ if and only if the value function $u^*$ is \lsc{} at $x$ and $\Phi^*(x)$ is nonempty.

\begin{theorem}
    [Necessary and sufficient conditions for {\lminsc{}}]
    For $x \in \Dom(\Phi)$, the following statements are equivalent:
    \begin{enumerate}
        \item The function $u : \Gr(\Phi) \subset X \times Y \to \bar \RR$ is \lminsc{} at $x$. \label{item:lminsc-def}
        \item There exists $y^* \in \Phi(x)$ such that for every net $\seq{(x_\alpha,y_\alpha)}_{\alpha \in I} \subset \Gr(\Phi)$ with $x_\alpha \to x$ satisfies $u(x,y^*) \leq \liminf_{\alpha \in I}       u(x_\alpha, y_\alpha)$. \label{item:lminsc-nets-b}
        \item The function $u^* : \Dom(\Phi) \subset X \to \bar\RR$ is \lsc{} at $x$, and $\Phi^*(x) \ne \emptyset$.\label{item:lminsc-lsc-min}
    \end{enumerate}
    \label{thm:lminsc}
\end{theorem}

\begin{proof}
\condition{\ref{item:lminsc-def}}{\ref{item:lminsc-lsc-min}} Let $y^* \in \Phi(x)$ be as in Definition~\ref{def:lminsc}, and let $\gamma \in \RR$. If $\gamma < u^*(x)$, fix $\gamma' \in (\gamma', u^*(x))$. Then there is a neighborhood $U(x)$ such that $\gamma' < u(x', y')$ for each $x' \in U(x)$ and $y' \in \Phi(x')$. Then $\gamma < \gamma' \leq u^*(x')$ for each $x' \in U(x)$, so $u^*$ is \lsc{} at $x$. Furthermore, since $\gamma < u(x, y^*)$ implies $\gamma < u(x,y)$ for each $y \in \Phi(x)$, it follows that $y^* \in \Phi(x)$.

\condition{\ref{item:lminsc-lsc-min}}{\ref{item:lminsc-nets-b}} Let $y^* \in \Phi^*(x)$. Then for any net $\seq{(x_\alpha, y_\alpha)}_{\alpha \in I} \subset \Gr(\Phi)$ with $x_\alpha \to x$, it follows that $u(x,y^*) = u^*(x) \leq \liminf_{\alpha \in I} u^*(x_\alpha) \leq \liminf_{\alpha \in I} u(x_\alpha,y_\alpha)$, as needed.

\condition{\ref{item:lminsc-nets-b}}{\ref{item:lminsc-def}} Suppose $u$ is not \lminsc{} at $x$. Then for each $y \in \Phi(x)$, there is an $\gamma \in \RR$ with $\gamma < u(x,y)$ such that for every neighborhood $U(x)$ there exist $x' \in U(x)$ and $y' \in \Phi(x')$ with $\gamma > u(x', y')$. Taking the directed set $I = \set{U \subset X : \text{$U$ is a neighborhood of $x$}}$ with $U \cleq V$ iff $V \subset U$, it follows that there exists a net $\seq{(x_\alpha,y_\alpha)}_{\alpha \in I} \subset \Gr(\Phi)$ such that $u(x,y) > \gamma > u(x_\alpha, y_\alpha)$ for each $\alpha \in I$. Thus, we proved that for each $y \in \Phi(x)$ there is a net $\seq{(x_\alpha,y_\alpha)}_{\alpha \in I} \subset \Gr(\Phi)$ with $x_\alpha \to x$ and $u(x,y) > \gamma \geq \sup_{\alpha \in I} u(x_\alpha,y_\alpha) \geq \liminf_{\alpha \in I} u(x_\alpha, y_\alpha)$. The direct implication follows from the contrapositive.
\end{proof}

Theorem~\ref{thm:lminsc} admits an additional necessary and sufficient condition when $u$ is globally \lminsc{}. 

\begin{corollary}
    Suppose $\Phi : X \to 2^Y$ is a strict multifunction. Then the following statements are equivalent:
    \begin{enumerate}
        \item The function $u : \Gr(\Phi) \subset X \times Y \to \bar \RR$ is \lminsc{} at each $x \in X$. \label{item:lminsc-def-2}
        \item The multifunction $\epi u^* : X \to 2^{\RR}$ has a closed graph, and the equality $P_X^\Phi \epi u = \epi u^*$ holds. \label{item:lminsc-epi}
        \item $u^*$ is \lsc{} on $X$, and $\Phi^* : X \to 2^Y$ is a strict multifunction. \label{item:lminsc-lsc-min-2}
    \end{enumerate}
    \label{cor:lminsc}
\end{corollary}

\begin{proof}
Conditions~\ref{item:lminsc-def-2} and~\ref{item:lminsc-lsc-min-2} are equivalent in view of Theorem~\ref{thm:lminsc}. Conditions~\ref{item:lminsc-lsc-min-2} and~\ref{item:lminsc-epi} are equivalent in view of two facts: $u^*$ is \lsc{} if and only if $\epi u^*$ has a closed graph, and $\Phi^*(x) \ne \emptyset$ if and only if there exists $y^* \in \Phi(x)$ such that $(\epi u)(x,y^*) = (\epi u^*)(x)$.
\end{proof}

As with Corollary~\ref{cor:linfsc}, there are two ways to apply Corollary~\ref{cor:lminsc} when $\Phi$ is not necessarily strict. First, because $\Phi$ is strict on $\Dom(\Phi)$, the results hold when $X$ is replaced with $\Dom(\Phi)$. Second, the original minimization problem can be reduced to a problem where $u(x,y)$ is set to $+\infty$ for $y \notin \Phi(x)$, and where the multifunction is redefined to equal $Y$ for each $x \in X$; see Corollary~\ref{cor:strict-phi-bar-lminsc}. We remark that Corollary~\ref{cor:lminsc}, which provides necessary and sufficient conditions to topological spaces $X$ and $Y$, generalizes Proposition 1.18 in~\citet{rockafellar2009variational}, which gives sufficient conditions for the epigraphical projection equality when $X$ and $Y$ are Euclidean spaces.

Now let us consider sufficient conditions for lower semicontinuity of $u^*$. The following definition was introduced in~\citet[Definition 1.3]{feinberg2014berge} and~\citet[Definition 1]{feinberg2015continuity} for Hausdorff topological spaces. Here we apply it to general topological spaces.

\begin{definition}
    [{\citet[Definition 1.3]{feinberg2014berge}}]
    Let $E \subseteq \Dom(\Phi)$.  The function $u : \Gr(\Phi) \subseteq X \times Y \to \bar \RR$ is called \emph{$\KK\NN$-inf-compact} on $\Gr_E(\Phi)$ if the following conditions hold:
    \begin{enumerate}
        \item $u$ is \lsc{} at each $(x,y) \in \Gr_E(\Phi)$;
        \item if a net $\seq{x_\alpha}_{\alpha \in I}$ with values in $\Dom(\Phi)$ converges to $x \in E$, then each net $\seq{y_\alpha}_{\alpha \in I}$ with $y_\alpha \in \Phi(x_\alpha)$ for each $\alpha \in I$ satisfying the condition that $\set{u(x_\alpha,y_\alpha)}_{\alpha \in I}$ is bounded above, has an accumulation point $y \in \Phi(x)$.
    \end{enumerate}
    The function $u$ is \emph{$\KK\NN$-sup-compact} on $\Gr_{E}(\Phi)$ if $-u$ is $\KK\NN$-inf-compact on $\Gr_{E}(\Phi)$.
    \label{def:kn-inf-compact}
\end{definition}

For $X = Y = \RR$ with $\Phi(x) \equiv Y$, an important example of a $\KK\NN$-inf-compact function is $u(x,y) = \abs{x-y}^p$ for $1 \leq p < +\infty$. Section~\ref{sec:inventory-control} describes $\KK\NN$-inf-compact functions (\ref{eq:inventory-objective}) for inventory control.

In general the nets in Definition~\ref{def:kn-inf-compact}(ii) may converge to multiple limits.  Let $E' \subset E \subset \Dom(\Phi)$. If $u$ is $\KK\NN$-inf-compact on $\Gr_{E}(\Phi)$, then $u$ is also $\KK\NN$-inf-compact on $\Gr_{E'}(\Phi)$. The converse can fail when, for example, $u$ is not \lsc{} on $E \setminus E'$. On the other hand, if $u$ is $\KK\NN$-inf-compact on $\Gr_{\set{x}}(\Phi)$ for each $x \in E$, then $u$ is $\KK\NN$-inf-compact on $\Gr_E(\Phi)$. In addition, $\KK\NN$-inf-compactness can be reformulated in the following manner: for each net $\seq{(x_\alpha,y_\alpha)_{\alpha \in I}} \subset \Gr(\Phi)$ with $x_\alpha \to x$, if $\sup_{\alpha \in I} u(x_\alpha, y_\alpha) \leq \lambda < +\infty$, then there exists a subnet $\seq{y_\beta}_{\beta \in J} \subset \seq{y_\alpha}_{\alpha \in I}$ and $y \in \Phi(x)$ such that $(x_\beta, y_\beta) \to (x,y)$ and $\liminf_{\beta \in J} u(x_\beta, y_\beta) \leq \lambda$. This formulation directly implies that $\KK\NN$-inf-compactness is a stronger property than lower inf-semicontinuity in view of Theorem~\ref{thm:linfsc}\ref{item:linfsc-nets}. The following theorem shows that a $\KK\NN$-inf-compact function is \lminsc{} (hence, \linfsc{}). An example of an \lminsc{} (hence, \linfsc{}) function that is not $\KK\NN$-inf-compact is provided in Example~\ref{ex:optimum-counterexample}.  A further characterization of $\KK\NN$-inf-compactness is developed in Theorem~\ref{thm:n-inf-forceful}\ref{item:n-inf-forceful}.

\begin{theorem}
    [Properties of $\KK\NN$-inf-compact functions]
    Let $E \subset \Dom(\Phi)$. Consider the objective function $u : \Gr(\Phi) \subseteq X \times Y \to \bar\RR$.
    \begin{enumerate}
        \item The function $u$ is $\KK\NN$-inf-compact on $\Gr_{E}(\Phi)$ if and only if, for each net $\seq{(x_\alpha,y_\alpha)}_{\alpha \in I}$ with values in $\Gr(\Phi)$ with $x_\alpha \to x \in E$ satisfying $\liminf_{\alpha \in I} u(x_\alpha, y_\alpha) < +\infty$, there exists a convergent subnet $\seq{y_{\beta}}_{\beta \in J}$ of the net $\seq{y_\alpha}_{\alpha \in I}$ with a limit $y \in \Phi(x)$, and each limit $y \in \Phi(x)$ of $\seq{y_\beta}_{\beta \in J}$ satisfies $u(x,y) \leq \liminf_{\alpha \in I} u(x_\alpha, y_\alpha)$. \label{item:n-inf-forceful}
        \item If $u$ is $\KK\NN$-inf-compact on $\Gr_E(\Phi)$, then $u$ is \lminsc{} at each $x \in E$.
    \end{enumerate}
    \label{thm:n-inf-forceful}
\end{theorem}

\begin{proof}
Suppose $u$ is $\KK\NN$-inf-compact on $\Gr_{E}(\Phi)$, and let $x \in E$. Let $\seq{(x_\alpha, y_\alpha)}_{\alpha \in I}$ be an arbitrary net with values in $\Gr(\Phi)$, with $x_\alpha \to x$, and satisfying $\liminf_{\alpha \in I} u(x_\alpha, y_\alpha) < +\infty$.  There exists a subnet $\seq{(x_\beta, y_\beta)}_{\beta \in J}$ of the net $\seq{(x_\alpha, y_\alpha)}_{\alpha \in I}$ with $\lim_{\beta \in J} u(x_\beta, y_\beta) = \liminf_{\alpha \in I} u(x_\alpha, y_\alpha)$. Fix $\beta \in J$ such that $\abs{u(x_{\beta'}, y_{\beta'}) - \lim_{\beta \in J} u(x_\beta, y_\beta)} < 1$ for each $\beta' \in J$ with $\beta \cleq \beta'$.  Let $J_\beta = \set{\beta' \in J : \beta \cleq \beta'}$. Since $u$ is $\KK\NN$-inf-compact, the subnet $\seq{(x_{\beta'},y_{\beta'})}_{\beta' \in J_\beta}$ satisfies $\sup_{\beta' \in J_\beta} u(x_{\beta'}, y_{\beta'}) < +\infty$, so there is a limit point $y \in \Phi(x)$ of the net $\seq{y_{\beta'}}_{\beta' \in J}$. Fix a subnet $\seq{y_\gamma}_{\gamma \in L}$ of the net $\seq{y_\beta}_{\beta \in J}$ which converges to this limit point. Since $u$ is \lsc{}, for each $y \in \Phi(x)$ with $y_\gamma \to y$, then $u(x,y') \leq \liminf_{\gamma \in L} u(x_\gamma, y_\gamma) = \liminf_{\alpha \in I} u(x_\alpha, y_\alpha)$, which is desired.

Conversely, to see that $u$ is \lsc{} on $\Gr_{E}(\Phi)$, let $x \in E$ and let $y \in \Phi(x)$. Suppose the net $\seq{(x_\alpha, y_\alpha)}_{\alpha \in I}$ with values in $\Gr(\Phi)$ converges $(x_\alpha, y_\alpha) \to (x,y)$. If $\liminf_{\alpha \in I} u(x_\alpha, y_\alpha) = + \infty$, then $u(x,y) \leq \liminf_{\alpha \in I} u(x_\alpha, y_\alpha)$ is immediate. Otherwise, there is a convergent subnet $\seq{y_\beta}_{\beta \in J}$ of the net $\seq{y_\alpha}_{\alpha \in I}$. Since $y$ is a limit of $\seq{y_\beta}_{\beta \in J}$, it follows that $u(x,y) \leq \liminf_{\alpha \in I} u(x_\alpha, y_\alpha)$, which shows that $u$ is \lsc{}. Next, let the net $\seq{x_\alpha}_{\alpha \in I}$ with values in $\Dom(\Phi)$ converge to $x$, and let the net $\seq{y_\alpha}_{\alpha \in I}$ with $y_\alpha \in \Phi(x_\alpha)$ for each $\alpha \in I$ and such that $\set{u(x_\alpha, y_\alpha)}_{\alpha \in I}$ is bounded above. Then $\liminf_{\alpha \in I} u(x_\alpha, y_\alpha) \leq \sup_{\alpha \in I} u(x_\alpha, y_\alpha) < +\infty$, so there exists a subnet $\seq{y_\beta}_{\beta \in J} \subseteq \seq{y_\alpha}_{\beta \in J}$ and $y \in \Phi(x)$ such that $y_\beta \to y$ and $u(x,y) \leq \liminf_{\alpha \in I} u(x_\alpha, y_\alpha)$. Therefore, $u$ is $\KK\NN$-inf-compact on $\Gr_{E}(\Phi)$.

Finally, suppose $u$ is $\KK\NN$-inf-compact on $\Gr_E(\Phi)$, and let $x \in E$. First, we observe that $\Phi^*(x) \ne \emptyset$. Indeed, if $u^*(x) = +\infty$, then $\Phi^*(x) = \Phi(x)$. Otherwise, let $\seq{y_n}_{n \in \NN} \subset \Phi(x)$ satisfy $u(x,y_n) \leq u^*(x) + n\inv$ for each $n$. Since $\sup_{n} u(x,y_n) < +\infty$, $\seq{y_n}_{n \in \NN}$ has an accumulation point $y^* \in \Phi(x)$ satisfying $u(x,y^*) \leq \liminf_{n \in \NN} u(x_n, y_n) = u^*(x)$, so $y^* \in \Phi^*(x)$. Next, let $\seq{(x_\alpha,y_\alpha)}_{\alpha \in I} \subset \Gr(\Phi)$ with $x_\alpha \to x$. If $\liminf_{\alpha \in I} u(x_\alpha,y_\alpha) = +\infty$, there is nothing to show. Otherwise, there exists an accumulation point $y \in \Phi(x)$ of the net $\seq{y_\alpha}_{\alpha \in I}$, which satisfies $u(x,y^*) = u^*(x) \leq u(x,y) \leq \liminf_{\alpha \in I} u(x_\alpha,y_\alpha)$, as needed.
\end{proof}

We consider a construction introduced in the discussion after Theorem 6~\citet{feinberg2015continuity}.  Let $x \in \Dom(\Phi)$ be fixed, and suppose there exist $y \in \Phi(x)$ and $\lambda \in \RR$ such that $u(x,y) < \lambda$.  Define the multifunction $\Phi_{\lambda,x} : X \to 2^Y$ by the equation
\begin{equation}
    \Phi_{\lambda,x}(z) =
    \begin{cases}
        \set{y \in \Phi(z) : u(z,y) \leq \lambda} & \text{if this set is nonempty,} \\
        \set{y \in \Phi(x) : u(x,y) \leq \lambda} & \text{otherwise,}
    \end{cases}
    \label{eq:modified-phi}
\end{equation}
and define $u_{\lambda,x} : \Gr(\Phi_{\lambda,x}) \subseteq X \times Y \to
\bar\RR$ by the equation
\begin{equation}
    u_{\lambda,x}(z,y) =
    \begin{cases}
        u(z,y) & \text{if $\set{y \in \Phi(z) : u(z,y) \leq \lambda} \ne \emptyset$,} \\
        u(x,y) & \text{otherwise.}
    \end{cases}
    \label{eq:modified-u}
\end{equation}
Let $u_{\lambda,x}^*$ be the value function for the minimization problem $(X, Y, \Phi_{\lambda,x}, u_{\lambda,x})$, and let $\Phi_{\lambda,x}^*$ be its solution multifunction.  The modified problem is related to the original problem $(X, Y, \Phi, u)$ through the next lemma. Both statements of Lemma~\ref{lem:relationship-of-minimization-problems} are obvious when $X$ and $Y$ are Hausdorff spaces and were used throughout~\citet{feinberg2015continuity} without proof.

\begin{lemma}
    Let $x \in X$, and suppose there exists $\lambda \in \RR$ with $u(x,y) < \lambda$ for some $y \in \Phi(x)$. The following statements hold:
    \begin{enumerate}
        \item The equalities $u_{\lambda,x}^*(z) = u^*(z)$ and $\Phi_{\lambda,x}^*(z) = \Phi^*(z)$ hold for all $z$ such that $u^*(z) < \lambda$. In particular, $u_{\lambda,x}^*(x) = u^*(x)$ and $\Phi_{\lambda,x}^*(x) = \Phi^*(x)$. Moreover, if $u_{\lambda,x}^*$ is \lsc{} at $x$, then $u^*$ is \lsc{} at $x$.
        \item If $u$ is $\KK\NN$-inf-compact on $\Gr_{\set{x}}(\Phi)$, then $u_{\lambda,x}$ is $\KK\NN$-inf-compact on $\Gr_{\set{x}}(\Phi_{\lambda,x})$.
    \end{enumerate}
    \label{lem:relationship-of-minimization-problems}
\end{lemma}

\begin{proof}
\begin{enumerate}[wide,topsep=0pt,parsep=0pt]
    \item When $u^*(z) < +\infty$, the equality $u_{\lambda,x}^*(z) = u^*(z)$ is evident from the definitions of $\Phi_{\lambda,x}$ and $u_{\lambda,x}$, and therefore the $\Phi_{\lambda,x}^*(z) = \Phi^*(z)$. In addition, $u_{\lambda,x}^*(z) \leq u^*(z)$ for each $z \in \Dom(\Phi)$. Therefore, if a net $\seq{x_\alpha}_{\alpha \in I}$ with values in $\Dom(\Phi)$ converges to $x$, the inequality $\liminf_{\alpha \in I} u^*(x_\alpha) \geq \liminf_{\alpha \in I} u_{\lambda,x}^*(x_\alpha) \geq u_{\lambda,x}^*(x) = u^*(x)$ implies that $u^*$ is \lsc{} at $x$.
    \item Let $\varphi : X \to X$ be defined by $\varphi(z) = z$ if $\mathcal{D}_u(\lambda; \Phi(z)) \ne \emptyset$ and $\varphi(z) = x$ otherwise. To see that $u_{\lambda,x}$ is \lsc{} at each $(x,y) \in \Gr_{\set{x}}(\Phi_{\lambda,x})$, let $\seq{(x_\alpha,y_\alpha)}_{\alpha \in I}$ be a net with values in $\Gr(\Phi_{\lambda,x})$ converging to $(x,y) \in \Gr_{\set{x}}(\Phi_{\lambda,x})$. We observe that the net $\seq{(\varphi(x_\alpha),y_\alpha)}_{\alpha \in I}$ is contained in $\Gr(\Phi)$ and converges to $(x,y) \in \Gr_{\set{x}}(\Phi)$, and moreso that $u_{\lambda,x}(x_\alpha,y_\alpha) = u(\varphi(x_\alpha), y_\alpha)$ for each $\alpha \in I$. The inequality $u_{\lambda,x}(x,y) = u(x,y) \leq \liminf_{\alpha \in I} u(\varphi(x_\alpha),y_\alpha) = \liminf_{\alpha \in I} u_{\lambda,x}(x_\alpha,y_\alpha)$ implies that $u_{\lambda,x}$ is \lsc{} at each $(x,y) \in \Gr_{\set{x}}(\Phi_{\lambda,x})$.  Next, let the net $\seq{(x_\alpha, y_{\alpha})}_{\alpha \in I}$ with values in $\Gr(\Phi_{\lambda,x})$ and with $x_\alpha \to x$ be given. Since this net is contained in $\Gr(\Phi_{\lambda,x})$, the inequality $\sup_{\alpha \in I} u_{\lambda,x}(x_\alpha,y_\alpha)  = \sup_{\alpha \in I} u(\varphi(x_\alpha), y_\alpha)\leq \lambda < +\infty$ holds. Therefore, the net $\seq{y_\alpha}_{\alpha \in I}$ has an accumulation point $y \in \Phi(x)$, and we can fix a subnet $\seq{(x_\beta,y_\beta)}_{\beta \in J}$ such that $(x_\beta, y_\beta) \to (x, y)$. Since $u$ is \lsc{} at $(x,y)$, there holds $u(x,y) \leq \liminf_{\beta \in I} u(\varphi(x_\beta), y_\beta) \leq \lambda$, which implies that $y \in \Phi_{\lambda,x}(x)$, and so $\seq{y_\alpha}_{\alpha \in I}$ has an accumulation point in $\Phi_{\lambda,x}(x)$. This verifies that $u_{\lambda,x}$ is $\KK\NN$-inf-compact on $\Gr_{\set{x}}(\Phi_{\lambda,x})$.  \qedhere
\end{enumerate}
\end{proof}

Although $\KK\NN$-inf-compactness is still a sufficient condition for $u^*$ to be \lsc{}, it can be further generalized by examining the modified minimization problem defined by the tuple $(X, Y, \Phi_{\lambda,x}, u_{\lambda,x})$  in (\ref{eq:modified-phi}) and (\ref{eq:modified-u}). The following theorem was proven by~\citet[Theorem 7]{feinberg2015continuity} when $X$ and $Y$ are Hausdorff spaces, and it provides a local sufficient condition for the function $u^*$ to be \lsc{} and the solution set $\Phi^*(x)$ to be nonempty and compact. Here we prove it for general topological spaces. It also generalizes Corollary 1 from~\citet{feinberg2015continuity} to general topological spaces.

\begin{theorem}
    [{cf.~\citet[Theorem 7]{feinberg2015continuity}}]
    Let $x \in \Dom(\Phi)$, and suppose there exists $\lambda \in \RR$ such that $u(x,y) < \lambda$ for some $y \in \Phi(x)$. For the minimization problem defined by the tuple $(X, Y, \Phi_{\lambda,x}, u_{\lambda,x})$ in (\ref{eq:modified-phi}) and (\ref{eq:modified-u}), if $u_{\lambda,x}$ is $\KK\NN$-inf-compact on $\Gr_{\set{x}}(\Phi_{\lambda,x})$, then $u^*$ is \lsc{} at $x$, and $\Phi^*(x)$ is nonempty and compact.
    \label{thm:general-lsc-value}
\end{theorem}

\begin{proof}
Since $u_{\lambda,x}$ is $\KK\NN$-inf-compact on $\Gr_{\set{x}}(\Phi_{\lambda,x})$, it follows from Theorem~\ref{thm:n-inf-forceful} that $u_{\lambda,x}$ is \lminsc{} at $x$, so $u_{\lambda,x}^*$ is \lsc{} at $x$. In view of Theorem~\ref{thm:lminsc} and Lemma~\ref{lem:relationship-of-minimization-problems}, it follows that $u^*$ is \lsc{} at $x$ and $\Phi^*(x) \ne \emptyset$.
To see that $\Phi^*(x)$ is compact, it suffices in view of Lemma~\ref{lem:relationship-of-minimization-problems} to show that $\Phi_{\lambda,x}^{*}(x)$ is compact.  To this end, let $\seq{y_\alpha}_{\alpha \in I}$ be a net with values in $\Phi_{\lambda,x}^*(x)$. Since $\sup_{\alpha \in I} u_{\lambda,x}(x,y_\alpha) \leq \lambda < +\infty$, there is an accumulation point $y \in \Phi_{\lambda,x}(x)$ of $\seq{y_\alpha}_{\alpha \in I}$.  According to Theorem~\ref{thm:n-inf-forceful}, this accumulation point satisfies $u_{\lambda,x}(x,y) \leq \liminf_{\alpha \in I}u_{\lambda,x}(x,y_\alpha) = u_{\lambda,x}^*(x)$, which implies that $y \in \Phi_{\lambda,x}^{*}$, verifying compactness.
\end{proof}

\subsection{Continuity of values and upper semicontinuity of solutions}
\label{subsection:main-results}

This section addresses continuity of the value function $u^*$ and upper semicontinuity of the solutions multifunction $\Phi^*$.  The following theorem is the characterization of continuity of $u^*$.

\begin{theorem}
    [Characterization of continuity of $u^*$; existence of solutions]
    The function $u^* : \Dom(\Phi) \subset X \to \bar \RR$ is continuous at $x \in \Dom(\Phi)$ if and only if $u$ is \FPTusc{} on $\Gr_{\set{x}}(\Phi)$ and \linfsc{} at $x$.  Furthermore, $u^*$ is continuous at $x$ and $\Phi^*(x) \ne \emptyset$ if and only if $u$ is \FPTusc{} on $\Gr_{\set{x}}(\Phi)$ and \lminsc{} at $x$.
    \label{thm:continuity-value}
\end{theorem}

\begin{proof}
This follows from combining Theorems~\ref{thm:general-usc-value},~\ref{thm:linfsc}, and~\ref{thm:lminsc}.
\end{proof}

Whereas Theorems~\ref{thm:general-usc-value} and~\ref{thm:lminsc} give necessary and sufficient conditions for $u^*$ to be \usc{} and \lsc{} with $\Phi^*(x) \ne \emptyset$, respectively, no necessary and sufficient conditions for $\Phi^*$ to be compact-valued and \usc{} are known.~\citet[Theorem 3]{tian1995transfer} provide necessary and sufficient conditions for $\Phi^*$ to be compact-valued and \usc{} when $\Phi$ is compact-valued with a closed graph. In the following theorem, we give general sufficient conditions for $\Phi^*$ to be \usc{} and compact-valued. This theorem generalizes ~\citet[Theorem 1]{tian1992maximum} by allowing possibly noncompact-valued $\Phi$.  It also generalizes ~\citet[Theorem 1.4]{feinberg2014berge}, by weakening the assumption that $\Phi$ is \lsc{} and $u$ is \usc{}.

\begin{theorem}
    [Local maximum theorem]
    Let $x \in \Dom(\Phi)$.  Suppose the function $u : \Gr(\Phi) \subset X \times Y \to \bar \RR$ is \FPTusc{} on $\Gr_{\set{x}}(\Phi)$ and $\KK\NN$-inf-compact on $\Gr_{\set{x}}(\Phi)$. Then the function $u^* : \Dom(\Phi) \subset X \to \bar \RR$ is continuous at $x$. If $u^*(x) = +\infty$, then $\Phi^*(x) = \Phi(x)$. Otherwise, $\Phi^*$ is \usc{} at $x$, and $\Phi^*(x)$ is nonempty and compact.
    \label{thm:main-berge-max}
\end{theorem}

We delay the proof until the end of this section. An immediate consequence of Theorem~\ref{thm:main-berge-max} is its global formulation.

\begin{corollary}
    [Maximum theorem]
    Suppose the function $u : \Gr(\Phi) \subset X \times Y \to \bar \RR$ is \FPTusc{} on $\Gr(\Phi)$ and $\KK\NN$-inf-compact on $\Gr(\Phi)$. Then the function $u^* : \Dom(\Phi) \subset X \to \bar \RR$ is continuous. For each $x \in \Dom(\Phi)$, if $u^*(x) = +\infty$, then $\Phi^*(x) = \Phi(x)$, and if $u^*(x) < +\infty$, then $\Phi^*$ is \usc{} at $x$, and $\Phi^*(x)$ is nonempty and compact.
    \label{cor:main-berge-max}
\end{corollary}

Theorem~\ref{thm:main-berge-max} is a direct application of the next result, which is the most general known form of Berge's maximum theorem. Theorem~\ref{thm:general-local-max} is a stronger result than Theorem 1 from~\citet{tian1992maximum} because it does not assume that the feasible sets $\Phi(x)$ are compact. It is also stronger than Theorem 7 from~\citet{feinberg2015continuity}, because it does not assume that $\Phi$ is \lsc{} or that $u$ is \usc{}, which together imply that $u$ is
\FPTusc{}.

\begin{theorem}
    [Generalized local maximum theorem]
    Let $x \in \Dom(\Phi)$.  Suppose $u$ is \FPTusc{} on $\Gr_{\set{x}}(\Phi)$, and suppose there exists $\lambda \in \RR$ such that $u(x,y) < \lambda$ for some $y \in \Phi(x)$. For the minimization problem defined by the tuple $(X, Y, \Phi_{\lambda,x}, u_{\lambda,x})$ in (\ref{eq:modified-phi}) and (\ref{eq:modified-u}), if $u_{\lambda,x}$ is $\KK\NN$-inf-compact on $\Gr_{\set{x}}(\Phi_{\lambda,x})$, Then $u^*$ is continuous at $x$, $\Phi^*(x)$ is nonempty and compact, and $\Phi^*$ is \usc{} at $x$.
    \label{thm:general-local-max}
\end{theorem}

\begin{proof}
The continuity of $u^*$ at $x$ follows from applying Theorems~\ref{thm:general-usc-value} and~\ref{thm:general-lsc-value}.  Since $u^*(x) < \lambda$, there exists a neighborhood $U(x)$ such that $u^*(z) < \lambda$ for each $z \in U(x)$. In view of Lemma~\ref{lem:relationship-of-minimization-problems}(i), it follows that $u^*(z) = u_{\lambda,x}^*(z)$ for each $z \in U(x)$, so $u_{\lambda,x}^*$ is also continuous at $x$.  Moreover, $\Phi^*(z) = \Phi_{\lambda,x}^*(z)$ for each $z \in U(x)$.  The compactness of $\Phi^*(x)$ follows from Theorem~\ref{thm:general-lsc-value}.

We first show that $\Phi_{\lambda,x}^*$ is \usc{} and compact-valued at $x$. According to~\citet[Theorem 17.16]{aliprantis1999infinite}, it suffices to show that for every net $\seq{(x_\alpha, y_\alpha)}_{\alpha \in I} \subset \Gr(\Phi_{\lambda,x}^*)$ such that $x_\alpha \to x$, there exists an accumulation point $y \in \Phi_{\lambda,x}^*(x)$ of the net $\seq{y_\alpha}_{\alpha \in I}$.  Indeed, since $\sup_{\alpha \in I} u_{\lambda,x}(x_\alpha, y_\alpha) \leq \lambda$, and since $u_{\lambda,x}$ is $\KK\NN$-inf-compact on $\Gr_{\set{x}}(\Phi_{\lambda,x})$, the accumulation point $y \in \Phi_{\lambda,x}(x)$ exists by Definition~\ref{def:kn-inf-compact}. In addition, since $u_{\lambda,x}^*$ is continuous at $x$, it follows that
\begin{equation*}
    u_{\lambda,x}^*(x) \leq u_{\lambda,x}(x,y) \leq \liminf_{\alpha \in I}
    u_{\lambda,x}(x_\alpha, y_\alpha) = \liminf_{\alpha \in I}
    u_{\lambda,x}^*(x_\alpha) = u_{\lambda,x}^*(x),
\end{equation*}
so $u_{\lambda,x}(x,y) = u_{\lambda,x}^*(x)$. Therefore, $y \in \Phi_{\lambda,x}^*(x)$, which shows that $\Phi_{\lambda,x}^*$ is \usc{} and compact-valued at $x$.

We next show that $\Phi^*$ is \usc{} at $x$. Let $G \subset Y$ be open such that $\Phi^*(x) \subset G$. Since $\Phi^*(x) = \Phi_{\lambda,x}^*(x)$, and since $\Phi_{\lambda,x}^*$ is \usc{} at $x$, there is a neighborhood $V(x)$ such that $\Phi_{\lambda,x}^{*}(z) \subset G$ for each $z \in V(x)$. The equality $\Phi^*(z) = \Phi_{\lambda,x}^*(z)$ for each $z \in U(x)$ implies that $\Phi^*(z) \subset G$ for all $z \in U(x) \cap V(x)$, which implies that $\Phi^*$ is \usc{} at $x$.
\end{proof}

\begin{proof}[Proof of Theorem~\ref{thm:main-berge-max}.]
Suppose first that $u^*(x) = +\infty$. Then $u^*$ is \lsc{} at $x$, and by Theorem~\ref{thm:general-usc-value} $u^*$ is also \usc{} at $x$; hence, continuous. Finally, $\Phi^*(x) = \Phi(x)$ holds, because $u(x,y) \geq u^*(x) = +\infty$ for each $y \in \Phi(x)$. Now, suppose $u^*(x) < +\infty$. Then we can fix $y \in \Phi(x)$ and $\lambda < + \infty$ such that $u(x,y) < \lambda$. In view of Lemma~\ref{lem:relationship-of-minimization-problems}, the minimization problem defined by the tuple $(X, Y, \Phi_{\lambda,x}, u_{\lambda,x})$ in~(\ref{eq:modified-phi}) and~(\ref{eq:modified-u}) satisfies the assumptions of Theorem~\ref{thm:general-local-max}, and the results follow.
\end{proof}

\subsection{Equivalence of constrained and unconstrained parametric optimization problems}
\label{subsection:formulations}

To close this section, we discuss the relationship between minimization problems where the feasible sets $\Phi(x)$ are specified and constrained and problems where the sets are unconstrained, i.e., $\Phi(x) \equiv Y$. This is the form of the minimization problem considered by~\citet[Theorems 1.17 and 7.41]{rockafellar2009variational} for Euclidean spaces. There is a natural transformation from the constrained problem to the unconstrained problem, namely, to consider the problem $(X,Y,\bar{\Phi}, \bar{u})$, where $\bar{\Phi}(x)=Y$ for all $x \in X$, and the function $\bar{u}$ is defined according to the equation
\begin{equation}
    \bar{u}(x,y) =
    \begin{cases}
        u(x,y) & \text{if $y \in \Phi(x)$,} \\
        +\infty & \text{otherwise.}
    \end{cases}
    \label{eq:ubar}
\end{equation}

Furthermore, every unconstrained problem has a natural constrained problem; i.e., from the problem $(X,Y, \bar{\Phi}, \bar{u})$ we can consider the constrained problem $(X,Y, \hat{\Phi}, \hat{u})$, where $\hat{\Phi}(x) = \set{y \in Y : \bar{u}(x,y) < +\infty}$ and $\hat{u} : \Gr(\hat{\Phi}) \subset X \times Y \to \RR$ is defined as $\hat{u}(x,y) = \bar{u}(x,y)$. Note that in general, $\Phi \ne \bar{\Phi}$ and $\Phi \ne \hat{\Phi}$.  Instead, in our framework, there may exist $y \in \Phi(x)$ with $u(x,y) = +\infty$. As Theorem~\ref{thm:extension-iff} shows, the value functions $u^*$, $\bar{u}^*$, and $\hat{u}^*$ of coincide. The following theorems show that all three problems are equivalent with respect to whether the objective functions are \linfsc{}, $\KK\NN$-inf-compact, or \FPTusc{}. The problems $(X, Y, \Phi, u)$ and $(X, Y, \bar \Phi, u)$ are also equivalent with respect to whether the objective functions are \lminsc{}. However, these are only equivalent to $(X, Y, \hat{\Phi}, \hat{u})$ in this respect when $\hat{\Phi}(x) \ne \emptyset$; otherwise, $\Phi^*(x)$ and $\bar{\Phi}^*(x)$ consist of points $y$ that give infinite values, while $\hat{\Phi}^*(x)$ is empty.

\begin{theorem}
    Consider the minimization problems defined by the tuples $(X,Y,\Phi,u)$, $(X, Y, \bar{\Phi}, \bar{u})$, and $(X, Y, \hat{\Phi}, \hat{u})$. Then equality of the value functions $u^* = \bar u^* = \hat u^*$ holds, and for each $x \in \Dom(\Phi)$, the following statements hold:
    \begin{enumerate}
        \item If one of $u, \bar u, \hat u$ is \FPTusc{} on $\Gr_{\set{x}}(\Phi)$, so are the other two.
        \item If one of $u, \bar u, \hat u$ is \linfsc{} at $x$, so are the other two.
        \item $u$ is \lminsc{} at $x$ if and only if $\bar u$ is \lminsc{} at $x$.
        \item When $\hat{\Phi}(x) \ne \emptyset$, then if one of $u, \bar u, \hat u$ is \lminsc{} at $x$, so are the other two.
    \end{enumerate}
    Therefore, if Theorem~\ref{thm:general-usc-value} or~\ref{thm:linfsc} holds for one of these problems, it holds for the other two. Theorems~\ref{thm:lminsc} and~\ref{thm:continuity-value} hold for the problem $(X, Y, \Phi, u)$ if and only if they hold for $(X, Y, \bar{\Phi}, \bar{u})$. On the other hand, if $\hat{\Phi}(x) \ne \emptyset$, then if Theorems~\ref{thm:lminsc} or~\ref{thm:continuity-value} hold for one of these problems, it holds for the other two.
    \label{thm:extension-iff}
\end{theorem}

\begin{proof}
The equalities $u^* = \bar{u}^* = \hat{u}^*$ hold, as $\Gr(\hat \Phi) \subset \Gr(\Phi) \subset \Gr(\bar \Phi)$ by construction. Then if there exists $y \in \hat \Phi(x)$, the value functions are finite and equal; otherwise, if $\hat \Phi(x) = \emptyset$, the value functions are all infinite.
\begin{enumerate}[%
        align=left,
        itemindent=2\parindent,
        labelwidth=\parindent,
        listparindent=\parindent,
        leftmargin=0pt,
        partopsep=0pt,
        topsep=0pt,
        parsep=0pt,
    ]
    \item  It suffices to consider only $\bar u$ and $\hat u$. In view of Theorem~\ref{thm:general-usc-value}, $\bar{u}$ is \FPTusc{} on $\Gr_{\set{x}}(\bar{\Phi})$ if and only if $\bar{u}^*$ is \usc{} at $x$, if and only if $\hat{u}$ is \FPTusc{} on $\Gr_{\set{x}}(\hat{\Phi})$.
    \item It suffices to consider only $\bar u$ and $\hat u$. In view of Theorem~\ref{thm:linfsc}, $\bar{u}$ is \linfsc{} at $x$ if and only if $\bar{u}^*$ is \lsc{} at $x$, if and only if $\hat{u}$ is \linfsc{} at $x$.
    \item In view of Theorem~\ref{thm:lminsc}, $\bar{u}$ is \lminsc{} at $x$ if and only if $\bar{u}^*$ is \lsc{} at $x$, if and only if $u$ is \lminsc{} at $x$.
    \item Suppose $\hat \Phi(x) \ne \emptyset$. Then there exists $y \in \Phi(x) \cap \bar{\Phi}(x) \cap \hat{\Phi}(x)$ such that $u(x,y) = \bar u(x,y) = \hat{u}(x,y) < +\infty$. In this case, $\bar{\Phi}^*(x) = \hat{\Phi}^*(x)$. In view of Theorem~\ref{thm:lminsc}, $\bar u$ is \lminsc{} at $x$ if and only if $\bar u^*$ is \lsc{} at $x$ and $\bar{\Phi}^*(x) = \hat{\Phi}^*(x) \ne \emptyset$, if and only if $\hat{u}$ is \lminsc{} at $x$.
\end{enumerate}
In view of the above arguments, Theorems~\ref{thm:general-usc-value},~\ref{thm:lminsc},~\ref{thm:linfsc}, and~\ref{thm:continuity-value} hold in the stated cases.
\end{proof}

Because $\bar{\Phi}$ is defined as a strict multifunction, Corollaries~\ref{cor:linfsc} and~\ref{cor:lminsc} always hold for the minimization problem $(X, Y, \bar{\Phi}, \bar{u})$. The next two corollaries summarize this observation.

\begin{corollary}
    For the minimization problem $(X, Y, \bar{\Phi}, \bar{u})$, the following statements are equivalent:
    \begin{enumerate}
        \item The function $\bar{u}$ is \linfsc{} at each $x \in X$.
        \item The multifunction $\epi \bar{u}^* : X \to 2^{\RR}$ has a closed graph.
        \item The function $\bar{u}^* : X \to \bar\RR$ is \lsc{} on $X$.
    \end{enumerate}
    \label{cor:strict-phi-bar-linfsc}
\end{corollary}

\begin{corollary}
    For the minimization problem $(X, Y, \bar{\Phi}, \bar{u})$, the following statements are equivalent:
    \begin{enumerate}
        \item The function $\bar{u}$ is \lminsc{} at each $x \in X$.
        \item The multifunction $\epi \bar{u}^* : X \to 2^{\RR}$ has a closed graph, and the equality $P_X^{\bar{\Phi}} \epi \bar{u} = \epi \bar{u}^*$ holds.
        \item $\bar{u}^*$ is \lsc{} on $X$, and $\bar{\Phi}^* : X \to 2^Y$ is a strict multifunction.
    \end{enumerate}
    \label{cor:strict-phi-bar-lminsc}
\end{corollary}

\begin{proof}[{Proof of Corollaries~\ref{cor:strict-phi-bar-linfsc} and~\ref{cor:strict-phi-bar-lminsc}.}]
The multifunction $\bar{\Phi}$ is strict, so the results follow from Corollaries~\ref{cor:linfsc} and~\ref{cor:lminsc}, respectively. \qedhere
\end{proof}

Finally, the next theorem shows that all three problems are equivalent under the conditions of Theorem~\ref{thm:general-local-max}.

\begin{theorem}
    Consider the minimization problems defined by the tuples $(X,Y,\Phi,u)$, $(X, Y, \bar{\Phi}, \bar{u})$, and $(X, Y, \hat{\Phi}, \hat{u})$. Then if the hypotheses of Theorem~\ref{thm:general-local-max} hold for one of these problems at the point $x \in X$, then they hold for the other two.
    \label{thm:extension-thm}
\end{theorem}

\begin{proof}
It suffices to consider only $(X, Y, \bar{\Phi}, \bar{u})$ and $(X, Y, \hat{\Phi}, \hat{u})$.  We first observe that the equality $\bar{u}^* = \hat{u}^*$ holds. This is clear when $\hat{\Phi}(x) \ne \emptyset$, but it holds in general because $\inf \emptyset = \inf \set{+\infty}$. In view of Theorem~\ref{thm:general-usc-value}, $\bar{u}$ is \FPTusc{} on $\Gr_{\set{x}}(\bar{\Phi})$ if and only if $\bar{u}^*$ is \usc{} at $x$, if and only if $\hat{u}$ is \FPTusc{} on $\Gr_{\set{x}}(\hat{\Phi})$.  If $\bar{u}$ is $\KK\NN$-inf-compact on $\Gr_{\set{x}}(\bar{\Phi})$, then we observe that $\hat{u}$ is \lsc{} on $\Gr_{\set{x}}(\hat{\Phi})$. There are two cases. If $\bar{u}(x,y) = +\infty$ for all $y \in \bar{\Phi}$, then $\hat{\Phi}(x) = \emptyset$, and $\hat{u}$ is $\KK\NN$-inf-compact on $\Gr_{\set{x}}(\hat{\Phi})$ vacuously. Otherwise, if $\seq{(x_\alpha, y_\alpha)}_{\alpha \in I} \subset \Gr(\hat{\Phi})$ is a net with $x_\alpha \to x$ and $\sup_{\alpha \in I} \hat{u}(x_\alpha, y_\alpha) < +\infty$, then $\seq{(x_\alpha,y_\alpha)}_{\alpha \in I} \subset \Gr(\bar{\Phi})$, so there exists an accumulation point $y \in \hat{\Phi}(x)$ of the net $\seq{y_\alpha}_{\alpha \in I}$. Moreover, since $\bar{u}$ is \lsc, $\bar{u}(x,y) \leq \sup_{\alpha \in I} \bar{u}(x_\alpha,y_\alpha) = \sup_{\alpha \in I} \hat{u}(x_\alpha, y_\alpha) < +\infty$, so $\bar{u}(x,y) = \hat{u}(x,y)$ and $y \in \hat{\Phi}$. Thus $\hat{u}$ is $\KK\NN$-inf-compact on $\Gr_{\set{x}}(\hat{\Phi})$.  Conversely, suppose $\hat{u}$ is $\KK\NN$-inf-compact on $\Gr_{\set{x}}(\hat{\Phi})$, and observe again that $\bar{u}$ is \lsc{} on $\Gr_{\set{x}}(\bar{\Phi})$. If $\seq{(x_\alpha, y_\alpha)}_{\alpha \in I} \subset \Gr(\bar{\Phi})$ with $x_\alpha \to x$ and $\sup_{\alpha \in I} \bar{u}(x_\alpha, y_\alpha) < +\infty$, then $y_\alpha \in \hat{\Phi}(x_\alpha)$ for each $\alpha \in I$, and so there is an accumulation point $y \in \hat{\Phi}(x) \subset \bar\Phi(x)$. Thus, $\bar{u}$ is $\KK\NN$-inf-compact on $\Gr_{\set{x}}(\hat{\Phi})$.
\end{proof}

\section{Continuity of minima in metric spaces} \label{sec:minimization-metric}

In this section we consider the minimization problem $(X, Y, \Phi, u)$ when $X$ and $Y$ are metric spaces. Although the metrics $d_X$ on $X$ and $d_Y$ on $Y$ are distinct, we will suppress the subscript when there is no ambiguity. Denote $B(x, \delta) = \set{x' \in X : d(x, x') < \delta}$ the ball of radius $\delta$ centered at $x \in X$.

Let $E \subset \Dom(\Phi)$. In metric spaces, the function $u$ is \FPTusc{} on $\Gr_{E}(\Phi)$ if and only if for each $x \in E$, $y \in \Phi(x)$, and $\gamma \in \RR$ with $u(x, y) < \gamma$, there exists $\delta > 0$ such that for each $x' \in B(x,\delta)$ there exists $y' \in \Phi(x')$ such that $u(x',y') < \gamma$. This is because the open balls $B(x,\delta)$ form a neighborhood base of each $x \in E$.  Feasible path transfer upper semicontinuity is closely related to epi-upper semicontinuity, defined in~\citet[Definition 7.39]{rockafellar2009variational} for Euclidean spaces. Generalized to metric spaces, epi-upper semicontinuity of $u$ is a geometric property that states that the epigraph multifunction $\epi u$ is \lsc{}.  Equivalently, $u$ is epi-upper semicontinuous at $x$ if and only if for each sequence $\seq{x_n}_{n = 1}^{\infty} \subset \Dom(\Phi)$ with $x_n \to x$ and $y \in \Phi(x)$, there is a sequence $y_n \to y$ such that $y_n \in \Phi(x_n)$ for each $n \in \NN$ and such that $\limsup_{n \to \infty} u(x_n, y_n) \leq u(x,y)$; see~\citet[pp. 270-272]{rockafellar2009variational} for a detailed explanation. The following theorem provides a sequential characterization of \FPTusc{} functions on metric spaces.

\begin{theorem}
    [Characterization of \FPTusc{} functions for metric spaces $X$ and $Y$]
    \label{thm:metric-fptusc}%
    Let $x \in X$. The objective function $u$ is \FPTusc{} on $\Gr_{\set{x}}(\Phi)$ if and only if, for each sequence $\seq{x_n}_{n=1}^\infty \subset \Dom(\Phi)$ with $x_n \to x$ there is a sequence $\seq{y_n}_{n = 1}^\infty$ such that $y_n \in \Phi(x_n)$ for each $n \in \NN$ and such that
    \begin{equation}\label{eq:fptusc-metric}%
        \limsup_{n \to \infty} u(x_n, y_n) \leq u(x,y).%
    \end{equation}%
\end{theorem}
The above theorem shows that \FPTusc{} of $u$, unlike epi-upper semicontinuity, does not require that each $y_n \to y$. Therefore, every epi-upper semicontinuous function is \FPTusc{}. As follows from Example~\ref{ex:fptusc-not-epi-usc}, below, an \FPTusc{} function may not be epi-upper semicontinuous.

\begin{proof}[Proof of Theorem~\ref{thm:metric-fptusc}.]
Suppose $u$ is \FPTusc{} on $\Gr_{\set{x}}(\Phi)$ and let $y \in \Phi(x)$. When $u(x,y) = +\infty$, the result follows immediately, so suppose $u(x,y) < +\infty$. Let $\seq{x_n}_{n=1}^\infty \subset \Dom(\Phi)$ converge to $x$. For each $k \in \NN$, the inequality $u(x,y) < u(x,y) + k\inv$ holds, and since $u$ is \FPTusc{}, we may fix $\delta_k > 0$ such that for each $x' \in B(x, \delta_k)$ there exists $y' \in \Phi(x')$ satisfying $u(x', y') < u(x,y) + k\inv$. Therefore, we may take the sequence $\seq{y_n^k}_{n=1}^{\infty}$ such that $y_n^k \in \Phi(x_n)$ for each $n \in \NN$ and such that $u(x_n, y_n^k) < u(x, y) + k\inv$. By taking the sequence $\seq{y_n^n}_{n=1}^{\infty}$ it follows that $\limsup_{n \to \infty} u(x_n, y_n^n) \leq \limsup_{n \to \infty} u(x,y) + n\inv = u(x,y)$, as needed.

Conversely, suppose $u$ is not \FPTusc{} on $\Gr_{\set{x}}(\Phi)$. Let $y \in \Phi(x)$ with $u(x,y) < \gamma < +\infty$ and let $\delta > 0$ be arbitrary. Then there exists $x' \in B(x, \delta)$ such that one has $\gamma \leq u(x_n, y')$ for all $y' \in \Phi(x')$. For $\delta_n = n\inv$ select $x_n$ so that for all $y' \in \Phi(x_n)$ one has $\gamma \leq u(x_n, y')$. Clearly, $x_n \to x$. On the other hand, if $\seq{y_n}_{n=1}^\infty$ is any sequence with $y_n \in \Phi(x_n)$ for each $n$, then it follows that $u(x_n, y_n) \geq \gamma$ for each $n$. Hence, $u(x,y) < \gamma \leq \limsup_{n \to \infty} u(x_n, y_n)$, for every such sequence $\seq{y_n}_{n=1}^{\infty}$, as needed.
\end{proof}

In a similar fashion to \FPTusc{} functions, a function $u$ is \lminsc{} at a point $x \in X$ if there exists $y^* \in \Phi(x)$ such that for each $\gamma \in \RR$ with $u(x,y^*) > \gamma$, there is a $\delta > 0$ such that for each $x' \in B(x,\delta)$ and $y' \in \Phi(x')$ one has $u(x',y') > \delta$. Furthermore, echoing Theorem~\ref{thm:lminsc}, $u$ is \lminsc{} at $x$ if and only if there is $y^* \in \Phi(x)$ such that for every sequence $\seq{(x_n,y_n)}_{n = 1}^{\infty} \subset \Gr(\Phi)$ with $x_n \to x$, the inequality $u(x,y^*) \leq \liminf_{n \to \infty} u(x_n,y_n)$ holds. Similarly, $u$ is \linfsc{} at $x$ if and only if for each $\gamma \in \RR$ with $\gamma < u(x,y)$ for every $y \in \Phi(x)$, there exists $\delta > 0$ such that $\gamma < u(x',y')$ for each $x' \in B(x,\delta)$ and $y' \in \Phi(x')$. Furthermore, in view of Theorem~\ref{thm:linfsc}, $u$ is \linfsc{} at $x$ if and only if for each $\lambda \in \RR$ and sequence $\seq{(x_n,y_n)}_{n \in \NN} \subset \Gr(\Phi)$ with $x_n \to x$, if $\liminf_{n \to \infty} u(x_n, y_n) < \gamma$, there exists $y \in \Phi(x)$ such that $u(x,y) < \gamma$.

The property of an \lminsc{} function can be interpreted as a generalization of epi-lower semicontinuous functions described by~\citet[Definition 7.39]{rockafellar2009variational}. Originally given on Euclidean spaces, epi-lower semicontinuity is a property that requires that $\Gr(\epi u)$ to be closed. Equivalently, $u$ is epi-lower semicontinuous if $u(x,y) \leq \liminf_{n \to \infty} u(x_n, y_n)$ for every sequence $(x_n, y_n) \to (x,y) \in \Gr(\Phi)$. In view of this characterization, the generalization from epi-lower semicontinuity to \lminsc{} can be seen as relaxing two conditions: only one $y^* \in \Phi(x)$ needs to be estimated, and the sequence $\seq{y_n}_{n \in \NN}$ does not need to converge to $y^*$.

The notion of $\KK$-inf-compactness is related to $\KK\NN$-inf-compactness and was formulated first in~\citet[{see Assumption \textbf{W}*}]{feinberg2012average} and~\cite{feinberg2013berge}.

\begin{definition}
    [$\KK$-inf-compactness]%
    \label{def:k-inf-compact-metric}%
    The function $u : \Gr(\Phi) \subset X \times Y \to \bar \RR$ is \emph{$\KK$-inf-compact} on $\Gr(\Phi)$ if for each $\lambda \in \RR$ and every nonempty, compact $K \subset \Dom(\Phi)$, the set $\DD_u(\lambda; \Gr_K(\Phi))$ is compact.
\end{definition}

The following definition renders a local version of $\KK\NN$-inf-compactness for metric spaces.

\begin{definition}
    [$\KK\NN$-inf-compactness in metric spaces]
    Let $E \subseteq \Dom(\Phi)$.  The function $u : \Gr_E(\Phi) \subseteq X \times Y \to \bar \RR$ is \emph{$\KK\NN$-inf-compact} on $\Gr_E(\Phi)$ if the following conditions hold:
        \begin{enumerate}
            \item $u$ is \lsc{} at each $(x,y) \in \Gr_E(\Phi)$;
            \item if a sequence $\seq{x_n}_{n=1}^{\infty}$ with values in $\Dom(\Phi)$ converges to $x \in E$, then each sequence $\seq{y_n}_{n=1}^\infty$ with $y_n \in \Phi(x_n)$ for each $n \in \NN$ satisfying the condition that $\set{u(x_n,y_n)}_{n=1}^{\infty}$ is bounded above, has a limit point $y \in \Phi(x)$.
        \end{enumerate}
    \label{def:kn-inf-compact-metric}
\end{definition}

Definitions~\ref{def:kn-inf-compact} and~\ref{def:kn-inf-compact-metric} are equivalent when $X$ and $Y$ are metric spaces.  This is true since a is $\KK$-inf-compact on $\Gr(\Phi)$ if it is $\KK\NN$-inf-compact on $\Gr(\Phi)$, and the converse implication holds under the following assumptions: (a) $X$ and $Y$ are Hausdorff topological spaces, and (b) $X$ is compactly generated, and metric spaces are compactly generated; see~\citet[Theorem 2.1]{feinberg2014berge} for details.

\citet{balder1992existence} gave a definition of a ``strongly coercive'' function on $B\subset X\times Y,$ and this definition was also employed by~\citet{dufour2018expected}. We use quotations because the same term was broadly used in the literature under a different meaning.  Any set $B\subset X\times Y$ determines the graph of a multifunction.  The function $u : B \subset X \times Y \to \bar{\RR}$ is \emph{``strongly coercive'' } if, for each sequence $\seq{(x_n,y_n)}_{n=1}^\infty$ with values in $B$ satisfying $x_n \to x \in E$ and $\liminf_{n \to \infty} u(x_n, y_n) < +\infty$, there exists a convergent subsequence $\seq{y_{n_k}}_{k=1}^\infty \subseteq \seq{y_n}_{n=1}^\infty$ with a limit $y \in \Phi(x)$ satisfying $u(x,y) \leq \liminf_{n \to \infty} u(x_n, y_n)$.  The following theorem, whose version for topological spaces is Theorem~\ref{thm:n-inf-forceful}, shows that the notions of ``strong coercivity'' and $\KK\NN$-inf-compactness on $\Gr_E(\Phi)$ coincide.

\begin{theorem}
    [Characterizations of $\KK\NN$-inf-compactness in metric spaces]
    Let $X$ and $Y$ be metric spaces, and let $E \subset \Dom(\Phi)$. Then $u$ is $\KK\NN$-inf-compact on $\Gr_E(\Phi)$, if and only if, for each sequence $\seq{(x_n,y_n)}_{n=1}^\infty$ with values in $\Gr(\Phi)$ satisfying $x_n \to x \in E$ and $\liminf_{n \to \infty} u(x_n, y_n) < +\infty$, there exists a convergent subsequence $\seq{y_{n_k}}_{k=1}^\infty \subseteq \seq{y_n}_{n=1}^\infty$ with a limit $y \in \Phi(x)$ satisfying $u(x,y) \leq \liminf_{n \to \infty} u(x_n, y_n)$.
    \label{thm:metric-characterization}
\end{theorem}

\begin{proof}
Suppose first that $u$ is $\KK\NN$-inf-compact on $\Gr_E(\Phi)$. Suppose the sequence $\seq{(x_n,y_n)}_{n = 1}^{\infty} \subset \Gr(\Phi)$ satisfies both $x_n \to x \in E$ and $\liminf_{n \to \infty} u(x_n,y_n) < +\infty$. Taking subsequences as needed, assume without loss of generality that $\lim_{n \to \infty} u(x_n,y_n) = \liminf_{n \to \infty} u(x_n, y_n)$ and that $\sup_{n \in \NN} u(x_n, y_n) < +\infty$. Then there exists a subsequence $\seq{y_{n_k}}_{k=1}^{\infty} \subset \seq{y_n}_{n=1}^{\infty}$ that converges to $y \in \Phi(x)$, and since $\lim_{k \to \infty} u(x_{n_k}, y_{n_k}) = \lim_{n \to \infty} u(x_n, y_n)$, it follows that $u(x,y) \leq \liminf_{n \to \infty} u(x_n, y_n)$.

Conversely, we first show that $u$ is \lsc{} at each $(x,y) \in \Gr_E(\Phi)$. Suppose the sequence $\seq{(x_n,y_n})_{n=1}^{\infty} \subset \Gr(\Phi)$ converges to $(x,y) \in \Gr_E(\Phi)$. If $\liminf_{n \to \infty} u(x_n, y_n) = +\infty$, there is nothing to show. Otherwise, $u(x,y) \leq \liminf_{n \to \infty} u(x_n, y_n)$ by assumption. Next, suppose $\seq{x_n}_{n=1}^{\infty} \subset \Dom(\Phi)$ with $x_n \to x \in E$, and let $y_n \in \Phi(x_n)$ ($n = 1,2,\dots$) be arbitrary such that $\sup_{n \in \NN} u(x_n, y_n) < +\infty$. Since $\liminf_{n \to \infty} u(x_n, y_n) \leq \sup_{n \to \infty} u(x_n, y_n) < +\infty$, there exists a subsequence $\seq{y_{n_k}}_{k =1}^{\infty} \subset \seq{y_n}_{n=1}^{\infty}$ such that $y_{n_k} \to y \in \Phi(x)$. This implies $\KK$-inf-compactness.
\end{proof}

\section{Counterexamples} \label{sec:examples}

Throughout this section, when the topological spaces $X$ and $Y$ are specified as subsets of $\RR$, we will always assume they are equipped with Euclidean topologies.  Theorem~\ref{thm:lminsc} gives necessary and sufficient conditions for the function $u^*$ to be \lsc{} and $\Phi^*$ to be nonempty; namely, that $u$ is \lminsc{}. Theorem~\ref{thm:linfsc} gives necessary and sufficient conditions for the function $u^*$ to be \lsc{}; i.e., that $u$ is \linfsc{}. While every \lminsc{} function is \linfsc{}, the reverse implication does not hold in the following example.

\begin{example}
    [A \linfsc{} function that is not \lminsc{}]
    Let $X = Y = [0,1]$, let $\Phi(x) = [0,1]$ for each $x \in X$, and consider the function $u(x,y) = \one_{\QQ}(x-y) + \one_{\set{0}}(y) + y$. For a fixed $x \in X$ and $\varepsilon > 0$, fix $z \in [x-\varepsilon, x) \setminus \QQ$, which exists by the density of the irrational numbers in $[x-\varepsilon ,x)$. Then with $y=x-z$, it follows that $0 < y \leq \varepsilon$ and $u(x,y) = y$. Thus, the value function $u^*(x) = 0$ for each $x \in X$, so $u$ is \linfsc{}. On the other hand, suppose for the sake of a contradiction that $y^* \in \Phi(x)$. Then $y^* \leq u(x,y^*)$, and $u(x,y^*) = u^*(x) = 0$, but $u(x,0) \geq 1$, a contradiction. Thus $\Phi^*(x) = \emptyset$ for each $x \in X$, so $u$ is not \lminsc{}.
    \label{ex:linfsc-not-lminsc}
\end{example}

The following example shows that $\Phi^*$ need not be compact-valued nor \usc{} under the condition that $u$ is \lminsc{}. Denote by $\QQ$ the set of rational numbers.

\begin{example}
    [Insufficiency of \lminsc{} for upper semicontinuity, compactness of solutions]
    Let $X = Y = [0,1]$, let $\Phi(x) = [0,1]$ for each $x \in X$, and consider the function
    \begin{equation*}
        u(x,y) = \one_{\QQ \times [0, \frac{1}{2}]}(x,y) + \one_{([0,1]\setminus\QQ) \times [\frac{1}{2},1]}(x,y)
    \end{equation*}
    Then $u^*(x) = 0$ for each $x \in X$, and $\Phi^*(x) = (\frac{1}{2},1]$ when $x \in \QQ$ and $\Phi^*(x) = [0,\frac{1}{2})$ otherwise. It is clear that $\Phi^*(x)$ is not compact for any $x \in X$. To see that $\Phi^*$ is not \usc{}, observe that for $x \in [0,1] \setminus \QQ$ there is a sequence of rational numbers $x_n \to x$; however, while $1 \in \Phi^*(x_n)$ for each $n$, one has $1 \notin \Phi^*(x)$. A similar argument holds for $x \in [0,1] \cap \QQ$. Therefore $\Phi^*$ is not upper semicontinuous at any $x \in X$.
\end{example}

In addition, while $u^*$ is \lsc{} if and only if $u$ is \lminsc{}, either inf-lower semicontinuity or min-lower semicontinuity neither implies nor is implied by classical \lsc{}.

\begin{example}
    [Independence of \lsc{} and \lminsc{} functions]
    Let $X = Y = [0,1]$, let $\Phi_1(x) = \set{1}$ if $x \leq \frac{1}{2}$ and $\Phi_1(x) = [0,1]$ otherwise, and let $u_1(x,y) = \one_{(\frac{1}{2},1]}(y)$. Then $u_1^*(x) = \one_{[0,\frac{1}{2}]}(x)$, which is not \lsc{}; hence, $u_1$ is \lsc{} but not \lminsc{}, and hence not \linfsc{}. On the other hand, letting $\Phi_2(x) = [0,1]$ for each $x \in X$, and letting $u_2(x,y) = \one_{\QQ \times \QQ}(x,y)$, it follows that $u_2^*(x) = 0$ for each $x \in X$, so $u_2$ is \lminsc{}, hence \linfsc{}, and not \lsc{}. In general, whether the objective function is \lsc{} is independent of the choice of $\Phi$, whereas \lminsc{} does depend on $\Phi$.
\end{example}

Theorem~\ref{thm:general-lsc-value} provides a sufficient local condition for $u^*$ to be \lsc{}, but the following example shows that the condition is not necessary. Hence, this property is stronger than \lminsc{}.

\begin{example}
    [Continuity of values and upper semicontinuity of solutions without $\KK\NN$-inf-compactness]
    Let $X = Y = [0,1]$, and consider
    \begin{equation*}
        \Phi(x) =
        \begin{cases}
            \set{1} & \text{if $x = 0$,} \\
             [0,1] & \text{otherwise}
        \end{cases}
    \end{equation*}
    with the function $u(x,y) = 1 - xy$. Observe that $u^*(x) = 1-x$ and $\Phi^*(x) = \set{1}$ for each $x \in X$, so in particular $u^*$ is continuous, $\Phi^*(x)$ is compact and nonempty for each $x \in X$, and $\Phi^*$ is \usc{}. We show that there does not exist $\lambda \in \RR$ such that $u_{\lambda, 0}$ is $\KK\NN$-inf-compact on $\Gr_{\set{0}}(\Phi_{\lambda, 0})$. To this end, it suffices to show that $u_{1, 0}$ is not $\KK\NN$-inf-compact on $\Gr_{\set{0}}(\Phi_{1,0})$. Since $\Phi_{1,0} = \Phi$, we observe that the sequence $\seq{(n\inv, \frac{1}{2})}_{n = 1}^{\infty} \subseteq \Gr(\Phi_{1,0})$ does not converge in $\Gr(\Phi_{1,0})$.
    \label{ex:optimum-counterexample}
\end{example}

The following example is relevant to Theorem~\ref{thm:metric-fptusc}.  It shows that a \FPTusc{} function may not be  epi-upper semicontinuous.

\begin{example}
    [A {\FPTusc{} function that is not epi-upper semicontinuous}]
    Let $X = Y = \RR$, and let $\Phi(x) = Y$ for each $x \in X$. Let
    \begin{equation*}
        u(x,y) = \one_{\set{x > 0, y > -1}}(x,y) + \one_{\set{x > 0, y > 1}}(x,y)
    \end{equation*}
    and observe that $u^*(x) = 0$ for each $x \in X$. According to Theorem~\ref{thm:general-usc-value}, the function $u$ is \FPTusc{}. For the sequence $x_n = n\inv$ ($n = 1,2,\dots$) converging to $0$ from the right, if $y_n \in \Phi(x_n)$ are chosen such that $y_n \to 1$, then the sequence $\seq{y_n}_{n = 1}^{\infty}$ is eventually in the set $(\frac{1}{2}, \frac{3}{2})$, and thus $u(x_n, y_n)$ is eventually the constant sequence $\seq{1}_{n = 1}^{\infty}$. Therefore, since $u(0,1)=0,$ there is no sequence $(n\inv, y_n) \to (0,1)$ satisfying $\limsup_{n \to \infty} u(n\inv, y_n) \leq u(0,1).$ So $u$ is not epi-upper semicontinuous.
    \label{ex:fptusc-not-epi-usc}
\end{example}

The following example  is due to~\citet{vasquez1995counterexample}. We recall that the function $u$ is \emph{inf-compact in $y$ on $\Gr(\Phi)$} if for every $\lambda \in \RR$ and $x \in X$, the set $\mathcal{D}_{u(x,\:\cdot\:)}(\lambda; \Gr(\Phi))$ is compact.  The example in~\cite{vasquez1995counterexample} shows that $u^*$ may not be \lsc{} for a continuous function $u$ that is inf-compact in $y$. We remark that $\KK\NN$-inf-compactness of $u$ on $\Gr(\Phi)$ is a slightly stronger assumption than the assumption that $u$ is \lsc{} and inf-compact in $y$ on $\Gr(\Phi)$, and $\KK\NN$-inf-compact functions are \lminsc{}, hence \linfsc{}; see Theorem~\ref{thm:n-inf-forceful}.

\begin{example}
    [Counterexample of~\citet{vasquez1995counterexample}]
    Let $X = \RR$, $Y = [0,+\infty)$, and let $\Phi(x)=Y$ for each $x \in X$. Let
    \begin{equation}
        u(x,y) =
        \begin{cases}
            1 + y & \text{if $x \leq 0$ or ($x > 0$, $0 \leq y \leq \frac{1}{2x}$)} \\
            2 + \frac{1}{x} - (2x+1)y & \text{if $x > 0$, $\frac{1}{2x} < y \leq \frac{1}{x}$} \\
            y - \frac{1}{x} & \text{if $x > 0$, $y > \frac{1}{x}$}
        \end{cases}
    \end{equation}
    The multifunction $\Phi$ is continuous, and the function $u$ is continuous and inf-compact in $y$. The value function $u^*(x) = \one_{\set{x \leq 0}}(x)$ is not \lsc{}, and the solutions multifunction $\Phi^*(x) = \set{0}$ when $x < 0$ and $\Phi^*(x) = \set{\frac{1}{x}}$ when $x \geq 0$. The function $u$ is not \linfsc{} at 0. To see why, let $x_n = \frac{1}{n}$, and let $y_n = n$ for each $n \in \NN$. Then $x_n \to 0$, and $\sup_{n \in \NN} u(x_n, y_n) = 0 < \frac{1}{2}$, but there does not exist $y \in \Phi(0)$ with $u(0, y) < \frac{1}{2}$. By Theorem~\ref{thm:linfsc}, $u$ is not \linfsc{}.
\end{example}

\section{Continuity of value functions for inventory control}
\label{sec:inventory-control}

The periodic review, single-item inventory control problem with ordering costs has a rich history of studies.  The problem was formulated by~\citet{arrow1960mathematical}, preceded by earlier studies initiated by~\citet{arrow1951optimal} and~\citet{dvoretzky1952inventory}, and presented in several monographs including~\citet{bensoussan2011dynamic}, \citet{zipkin2000foundations},~\citet[Chapter 11]{simchi-levi2005logic},~\citet{porteus2002foundations}, and~\citet{heyman2004stochastic}. It has been known for long time that the value function for a finite-horizon problem with backorders, unlimited orders, and unlimited storage capacity is continuous; see e.g., \citet{heyman2004stochastic} or~\citet[Theorem 8.3.4]{simchi-levi2005logic}, in spite of the discontinuity of the one-step costs caused by the presence of ordering costs. There was a conjecture that this continuity is due to continuity properties of $K$-convex functions, introduced by~\citet{arrow1960mathematical}; see~\citet[Lemma 7.1(e)]{heyman2004stochastic}.  This conjecture was disproved in~\citet[Example 5.1]{feinberg2017structure}.  In this subsection we show that this continuity follows from Berge's maximum theorem for problems with discontinuous objective functions.  We consider problems with bounded or unlimited order sizes and with bounded or unlimited storage capacity.  For simplicity we deal with finite-horizon problems. For the infinite-horizon version of the problem with backorders, unlimited order sizes, and unbounded storage capacity, the value function is also continuous; see~\citet[Theorem 5.3]{feinberg2017structure} for details.

\begin{figure}[t]
    \centering
    \includegraphics[width=\textwidth]{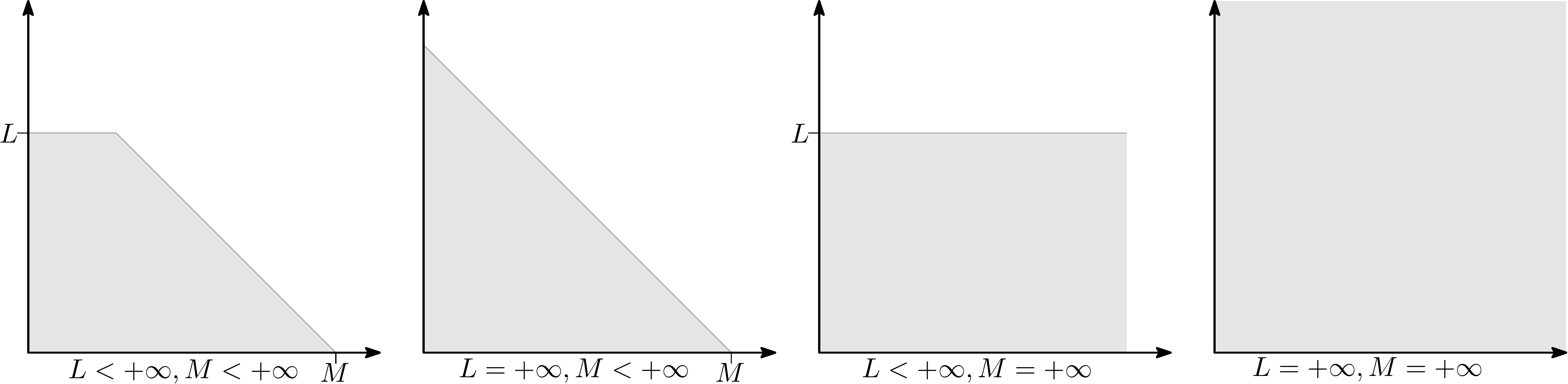}
    \caption{%
        Graph of feasible actions under four variants of the inventory control model.  Far left: finite storage capacity, limited orders. Center left: finite storage capacity, unlimited orders. Center right: infinite storage capacity, limited orders. Far right: infinite storage capacity, unlimited orders.
    }
    \label{fig:inventory-models}
\end{figure}

We consider the problem with possibly finite order size $L > 0$, possibly finite storage capacity $M > 0$, and possibly with lost sales. In the case of lost sales, we let $S := [0,+\infty)$, and, in the case of backorders, we let $S := \RR$. Without loss of generality, we assume $L \geq 1$. The choice of the state space $X$ depends on the value of $M$, and the choice of the action space $Y$ depends on the value of $L$.  If $M < + \infty$, then $X := S \cap (-\infty, M]$; otherwise, $X := S$.  Similarly, $Y := [0,L]$ when $L <  +\infty$ and $Y := [0,+\infty)$ otherwise.  The multifunction $\Phi : X \to 2^Y$ are specified differently depending on $L$ and $M$, as outlined below:
\begin{enumerate}
    \item Bounded orders, bounded storage capacity ($L, M < + \infty$): $\Phi(x) = [0, L \land (M-x)]$.
    \item Bounded orders, unlimited storage capacity ($L < +\infty$, $M = \infty$): $\Phi(x) = [0, L]$.
    \item Unlimited orders, bounded storage capacity ($L = +\infty$, $M < \infty$):  $\Phi(x) = [0, M-x]$.
    \item Unlimited orders, unlimited storage capacity ($L, M = +\infty$): $\Phi(x) = [0, +\infty)$.
\end{enumerate}
The following analysis will cover all four cases (with or without lost sales) without loss of generality.  The graphs of $\Phi$ in each of these cases are given in Figure~\ref{fig:inventory-models}. The analysis of each of these cases leads to the study of minimization problems over Euclidean spaces. However, different inventory control models may provide the controller only incomplete information of the inventory level; see, e.g.,~\citet{bai2021average},~\citet{bensoussan2011filtering}, and~\citet{feinberg2016partially}. Such models lead to minimization problems over infinite dimensional spaces. The analysis of these more complex models is broadly similar to the following finite-dimensional analysis, but in this paper we consider the latter for simplicity.

Let the function $T(x) = x^+$ when there are lost sales, and let $T(x) = x$ when there are backorders. We observe that $T$ is convex and nondecreasing. The state evolution is specified by the transition dynamics
\begin{equation*}
    x_{t+1} = T(x_t + y_t - D_{t+1}),
\end{equation*}
where $\seq{D_{t}}_{t = 1}^{\infty}$ are i.i.d.\ nonnegative demand shocks corresponding to the measure $P$.

For the finite-horizon, expected total discounted rewards criterion with discount factor $\alpha \in (0,1]$, we consider for $t = 0, 1, 2, \dots$ the minimization problems $(X, Y, \Phi, u_t)$, where the objective functions $u_t$ are given by $u_0(x,y) = 0$ and
\begin{equation}
    u_{t+1}(x,y) = c_1 \one_{\set{y > 0}}(y) + c_2 y +
    \EE h(T(x+y-D)) + \alpha
    \EE u_t^*(T(x+y-D))
    \label{eq:inventory-objective}
\end{equation}
Here, $c_1, c_2 > 0$, $h : \RR \to [0, +\infty)$ is a convex function satisfying $h(0) = 0$ and $\lim_{\abs{z} \to +\infty} h(z) = +\infty$, $u_t^*$ is the value function corresponding to $(X, Y, \Phi, u_t)$ (note that $u_0^*(x) = 0$ from the definition of $u_0(x,y)$), and $D$ is a representative random variable distributed identically to $\seq{D_t}_{t=1}^{\infty}$. We assume that $\EE h(z-D) < +\infty$ for each $z \in X$, which implies that $\EE h(z-D)$ is a continuous and convex function of $z$ on $X$.

If $L=M=+\infty$ and backorders are allowed, we deal with the classic well-studied case introduce by~\citet{arrow1960mathematical}  and clarified by~\citet{zabel1962note}.  Under certain conditions, there exist optimal $(s,S)$-policies which at time $n=0,\ldots,t-1$ place an order of the size $(S-x_t)$ if $x_t\le s$ and do not order otherwise.  In general, an optimal policy can have a different structure, and the structures of all optimal policies are described in~\citet{feinberg2017structure}.  Here we show that continuity of the values functions $u_t^*$ follows from Berge's maximum theorem for noncompact decision sets described in this paper. The continuity holds with or without backorders, for bounded or unlimited storage capacity, and with bounded or unlimited order sizes. The following two auxiliary lemmas are used to prove the main result, Theorem~\ref{thm:inventory-control}. The first gives sufficient conditions for the functional $U_t(z) = \EE u_t^*(T(z-D))$ to be continuous.

\begin{lemma}
    If a continuous function $f : X \to [0,+\infty)$ is bounded above by a continuous, convex function $g : X \to [0,+\infty)$, and $\EE g(T(z-D)) < +\infty$ for each $z \in X$, then the functional $F(z) = \EE f(T(z-D))$ is continuous for all $z \in X$.
    \label{lem:dominated-sequence}
\end{lemma}

\begin{proof}
We observe that $g \circ T$ is convex, continuous, and finite because $T$ is convex, continuous, finite, and nondecreasing. Furthermore, since $g \circ T$ is convex, the functional $G(z) = \EE g(T(z-D))$ is convex and thus continuous on $X$. In addition, $f \circ T \leq g \circ T$.  Let $z \in X$, and suppose $\seq{z_n}_{n =1}^\infty \subset X$ converges to $z$.  Then $f(T(z_n - D)) \to f(T(z-D))$ and $g(T(z_n-D)) \to g(T(z_n-D))$ almost surely, since $f$, $g$, and $T$ are continuous. Since $0 \leq f(T(z_n - D)) \leq g(T(z_n - D))$ and $0 \leq f(T(z-D)) \leq g(T(z-D))$, and since $\EE g(T(z_n - D)) \to \EE g(T(z-D))$, we apply the dominated convergence theorem to conclude that $\EE f(T(z_n - D)) \to \EE f(T(z-D))$. This shows that $F$ is continuous.
\end{proof}

The next technical lemma describes the construction of a feasible path that passes through fixed state-action pairs. This path is useful in the main theorem.

\begin{lemma}
    For each state-action pair $(x,y) \in \Gr(\Phi)$ and any sequence of states $\seq{x_n}_{n = 1}^{\infty} \subseteq X$ which converges to $x$, there exists a sequence of orders $\seq{y_n}_{n=1}^{\infty}$ such that $y_n \in \Phi(x_n)$ for each $n \in \NN$ which converges to $y$. In the case $y = 0$, the sequence $y_n = 0$ can be chosen for each $n \in \NN$.
    \label{lem:feasible-sequence}
\end{lemma}

\begin{proof}
Consider the function $\sigma_y : X \to Y$ defined by $\sigma_y(z) = y \one_{(-\infty, x)}(z) + (y - z + x) \one_{[x,x+y]}(z)$. We observe that $\sigma_y \in \Sigma(\Phi)$ and that $\sigma_y$ is continuous. Let $\seq{x_n}_{n = 1}^{\infty} \subseteq X$ be any sequence of states with $x_n \to x$. Then by continuity of $\sigma_y$, it follows that $\sigma_y(x_n) \to \sigma_y(x)$. Since $\sigma_y(x) = y$, $\seq{\sigma_y(x_n)}_{n=1}^{\infty}$ therefore suffices as a convergent sequence of feasible actions. Finally, observe that $\sigma_y(z) = 0$ whenever $y = 0$. This completes the proof.
\end{proof}

\begin{theorem}
    In the inventory control problem with possibly finite order size and possibly finite storage capacity, for each $t = 0, 1, 2, \dots$, the finite horizon expected total discounted value function $u^*_{t} : X \to \RR$ is continuous and satisfies the optimality equations
    \begin{align}
        u^*_0(x) &= 0, \label{eq:optimality-0} \\
        u^*_{t+1}(x) &= \min_{y \in \Phi(x)} \set{
            c_1 \one_{\set{y > 0}}(y) + c_2 y + \EE h(T(x+y-D))
            + \alpha \EE u_{t}^*(T(x+y-D))
        }
        \label{eq:optimality-equation}
    \end{align}
    Furthermore, there exists an optimal policy $\pi^*_t = (\varphi_0^*, \dots, \varphi_t^*)$ that achieves $u^*_t$.
    \label{thm:inventory-control}
\end{theorem}

The standard approach to proving continuity of the functions $u_t^*$ is by using the existence of optimal $(s,S)$ policies. This is the approach taken, for example, in~\citet{feinberg2016optimality, feinberg2017structure, feinberg2017optimality}. For problems with limited order sizes and capacities, optimal $(s,S)$ policies may not exist. The advantage of Theorem~\ref{thm:inventory-control} is that it avoids the analysis of optimal policies by employing Theorem~\ref{thm:main-berge-max} directly.

\begin{proof}[Proof of Theorem~\ref{thm:inventory-control}]
We will also prove that $U_t(z) = \EE u^*_t(T(z-D))$ is continuous for each $t = 0, 1, \dots$, which will follow from Lemma~\ref{lem:dominated-sequence} if we can show that $u_t^*(x) \leq g_t(x)$, where $g_t : X \to \RR$ is a nonnegative, continuous convex function, for each $t$.  We proceed by induction on $t$.  Since $u_0(x,y) \equiv 0$, then $u_0^*(x) \equiv 0$ for all $x \in X$, so $f_0(x) \equiv 0$ suffices.  Furthermore, any $y \in \Phi(x)$ achieves $u_0(x,y) = u_0^*(x)$.  Next, assume the result holds for $t$.

To show that $u_{t+1}^*$ is continuous, we first show that $u_{t+1}$ is $\KK$-inf-compact. To this end, we observe that $u_{t+1}(x,y)$ is \lsc{}. Let $K \subset X$ be nonempty and compact, and let $\lambda \in \RR$.  It follows that $\DD_{u_{t+1}}(\lambda; \Gr_K(\Phi))$ is a closed subset of the compact set $K \times [0, c_1\inv \lambda]$; hence, compact. According to Definition~\ref{def:k-inf-compact-metric}, $u_{t+1}$ is $\KK$-inf-compact.  We next show that $u_{t+1}$ is \FPTusc{} on $\Gr(\Phi)$. Given $(x,y) \in \Gr(\Phi)$, let $\seq{x_n}_{n=1}^{\infty} \subset X$ be arbitrary such that $x_n \to x$. According to Lemma~\ref{lem:feasible-sequence}, there exists a sequence $y_n \to y$ with $y_n \in \Phi(x_n)$ for each $n$. If $y > 0$, then $u_{t+1}$ is locally continuous at $(x,y)$, so $u_{t+1}(x_n, y_n) \to u_{t+1}(x, y)$.  If $y = 0$, then $u_{t+1}(\:\cdot\:, 0)$ is continuous at $x$. According to Lemma~\ref{lem:feasible-sequence}, we can select the feasible sequence $y_n \equiv 0$, in which case $u_{t+1}(x_n, 0) \to u_{t+1}(x, 0)$. Thus, Theorem~\ref{thm:metric-fptusc} implies that $u_{t+1}$ is \FPTusc{}. Then $u_{t+1}^*$ is continuous by applying Theorem~\ref{thm:main-berge-max} to the tuple $(X, Y, \Phi, u_{t+1})$. Let $\Phi^*_{t+1}(x) := \set{y \in \Phi(x) : u_{t+1}(x,y) = u_{t+1}^*(x)}$. Then Theorem~\ref{thm:main-berge-max} implies that $\Phi^{*}_{t+1}(x) \ne \emptyset$ for each $x \in X$, so we can fix a selector $\varphi^*_{t+1} \in \Sigma(\Phi^*_{t+1})$ that attains $u_{t+1}^*$. The optimality equation is satisfied by definition of the minimization problem in Equation~(\ref{eq:value-fn}) in which the infimum is achieved.

To see that $u^*_{t+1}(x) \leq g_{t+1}(x)$ for all $x \in X$, consider the policy $\varphi(x) \equiv 0$ that never orders. (This policy is available to the controller in each of the inventory control models.) From this policy, it follows that
\begin{equation*}
    u^*_{t+1}(x)  \leq u_{t+1}(x,\varphi(x))
    = \EE h(T(x - D)) + \alpha \EE u^*_t(T(x-D))
    \leq \EE h(T(x-D)) + \alpha \EE g_t(T(x-D)).
\end{equation*}
Let $g_{t+1}(x) := \EE h(T(x-D)) + \alpha \EE g_t(T(x-D))$, and observe that $g_{t+1}$ is the sum of two convex functions. Since $u^*_{t+1} \leq g_{t+1}$, the function $U_{t+1}(z) =  \EE u_{t+1}^{*}(T(z-D))$ is continuous in view of Lemma~\ref{lem:dominated-sequence}.

Therefore, the value functions $u^*_t$ are continuous for each $t = 0,1,2, \dots$, satisfy the optimality equations (\ref{eq:optimality-0}) and (\ref{eq:optimality-equation}), and the policy $\pi^*_t = (\varphi_0^*, \dots, \varphi_t^*)$ is optimal for each $t$.
\end{proof}

\section{Robust optimization} \label{sec:robust-optimization}

In this section we consider the \emph{minimax problem} $(X, A, B, f, \Phi_A, \Phi_B)$ defined by the equation
\begin{equation}
    f^*(x) = \adjustlimits \inf_{a \in \Phi_A(x)} \sup_{b \in \Phi_B(x,a)} f(x, a, b),
    \qquad x \in X,
    \label{eq:minimax-def}
\end{equation}
for topological spaces $X$, $A$, and $B$, a strict multifunction $\Phi_A : X \to 2^A$, a strict multifunction $\Phi_B : \Gr(\Phi_A) \subseteq X \times A \to 2^B$, and objective function $f : \Gr(\Phi_B) \subseteq X \times A \times B \to \bar{\RR}$. The function $f^*$ is called the \emph{minimax function}.

A minimax problem has the following game-theoretic interpretation. There are two players, \textsf{A} and \textsf{B}. Given the state $x \in X$ of the game, player \textsf{A} selects a feasible action $a \in \Phi_A(x)$. Then, after observing both the state $x$ and player \textsf{A}'s action $a$, player \textsf{B} selects a feasible action $b \in \Phi_B(x,a)$. The game ends with player \textsf{A} paying player \textsf{B} the quantity $f(x,a,b)$. Such games are called \emph{sequential} in the literature, and are a special case of \emph{perfect information} games; see, e.g.,~\citet[Theorem 1]{jaskiewicz2011stochastic} and~\cite[Section 1]{jaskiewicz2018zero}. Minimax problems have important applications to robust optimization. In this setting, the state information is partially uncertain, where $x$ represents certain data and $b$ represents uncertain data. The decision maker chooses an action $a \in \Phi_A(x)$ that ensures an upper bound for the cost $f(x,a,b)$ under all possible realizations of the uncertain data $b$. This is equivalent to the assumption that $b$ is selected by an adversary in view of the information $(x,a)$ within a specified uncertainty set $\Phi_B(x,a)$. This interpretation is employed in the study of robust dynamic programming or robust Markov decision processes, where the controller selects actions in the context of uncertain transition probabilities. Whereas in conventional dynamic programming, the Bellman equation for the value function resembles the form~(\ref{eq:value-fn}), the robust equivalent has the form~(\ref{eq:minimax-def}); see, e.g.,~\citet{iyengar2005robust},~\citet{nilim2005robust}, and~\citet{wiesemann2013robust} for analysis of models of this type. We remark that for robust Markov decision processes the function $f$ has a special structure.  We do not assume this structure here, and, in principle, our continuity results are applicable to other robust optimization problems.

The highest payoff player \textsf{B} can receive is the \emph{worst-loss function} $f^{\#} : \Gr(\Phi_A) \subseteq X \times A \to \bar{\RR}$ given by
\begin{equation}
    f^{\#}(x, a) = \sup_{b \in \Phi_B(x,a)} f(x,a,b),
    \qquad x \in X, \, a \in \Phi_A(x)
    \label{eq:worst-loss-sup}
\end{equation}
The smallest payment player \textsf{A} can make, assuming that player \textsf{B} maximizes her payoff, is the function $f^*$ defined in~(\ref{eq:minimax-def}), which also satisfies
\begin{equation}
    f^*(x) = \inf_{a \in \Phi_A(x)} f^{\#}(x,a).
    \label{eq:minimax-inf}
\end{equation}
Associated with $f^{\#}$ and $f^{*}$ are the \emph{optimal strategies} of players \textsf{A} and \textsf{B}, respectively:
\begin{align*}
    \Phi_B^\#(x,a) &= \set{b \in \Phi_B(x,a) : f^{\#}(x,a) = f(x,a,b)} \\
    \Phi_A^*(x) &= \set{a \in \Phi_A(x) : f^{*}(x) = f^{\#}(x,a)}.
\end{align*}
\citet{feinberg2017continuity, feinberg2018example} provide general sufficient conditions for the functions $f^*$ and $f^{\#}$ to be continuous and existence, compactness, and upper semi-continuity of optimal actions $\Phi_A^*$ and $\Phi_B^{\#}$ for players \textsf{A} and \textsf{B}, respectively. This section is a generalization of the prior work in three directions: it presents general necessary and sufficient conditions for continuity of $f^*$ and existence of optimal actions, the results are given in local form, and they are given for general topological spaces.

In the following results, we shall use the multifunctions $\Phi_A^{A \leftrightarrow B} : X \to 2^B$ and $\Phi_B^{A \leftrightarrow B} : \Gr(\Phi_A^{A \leftrightarrow B}) \subseteq X \times B \to 2^A$ defined by
\begin{align*}
    \Phi_A^{A \leftrightarrow B}(x) &= \set{b \in B : \text{$b \in \Phi_B(x,a)$
    for some $a \in \Phi_A(x)$}} \\
    \Phi_B^{A \leftrightarrow B}(x,b) &= \set{a \in A : b \in \Phi_B(x,a)}.
\end{align*}
The multifunction $\Phi_B^{A \leftrightarrow B}$ is closely related to the classical lower inverse of $\Phi_B$; see~\citet[Chapter 2, \S 3]{berge1963topological}. It can be understood as a kind of partial lower inverse. Similarly, we let $f^{A \leftrightarrow B} : X \times B \times A \to \bar{\RR}$ be defined by $f^{A \leftrightarrow B}(x,b,a) = f(x,a,b)$.

In the context of minimax problems, there is an additional semicontinuity assumption known as $A$-lower semicontinuity on $\Phi_B$ that provides the groundwork for the main results. This property was introduced in~\citet[Definition 4]{feinberg2017continuity} for the case that $X$, $A$, and $B$ are metric spaces.  The following definition gives the property of $A$-lower semicontinuity for arbitrary topological spaces.

\begin{definition}
    [$A$-lower semicontinuity]
    Let $E \subseteq \Dom(\Phi_A)$. The multifunction $\Phi_B : \Gr(\Phi_A) \subseteq X \times A \to 2^B$ is \emph{$A$-lower semicontinuous on $\Gr_E(\Phi_B)$} if for any net $\seq{(x_\alpha, a_\alpha)}_{\alpha \in I} \subseteq \Gr(\Phi_A)$ with $x_\alpha \to x$, $a \in A$ with $(x,a) \in E$, and $b \in \Phi_B(x,a)$ is arbitrary, there exists a net $\seq{b_\alpha}_{\alpha \in I}$ such that $b_\alpha \in \Phi_B(x_\alpha, a_\alpha)$ for each $\alpha$, and such that $b$ is a limit point of $\seq{b_\alpha}_{\alpha \in I}$.
    \label{def:a-lsc}
\end{definition}

\begin{remark}
    $A$-lower semicontinuity is an essential assumption for theorems concerning upper semicontinuity of minimax solutions; see, e.g.,~\citet{feinberg2018example}.
    If $\Phi_B$ is $A$-lower semicontinuous, then it is also \lsc{}. This is because, if $\seq{(x_\alpha, a_\alpha)}_{\alpha \in I} \subseteq \Gr(\Phi_A)$ and $(x_\alpha, a_\alpha) \to (x,a) \in \Gr(\Phi_A)$, then for each $b \in \Phi_B(x,a)$ there exists a net $\seq{b_\alpha}_{\alpha \in I}$ with $b_\alpha \in \Phi_B(x_\alpha, a_\alpha)$ for each $\alpha \in I$ having $b$ as a limit point. According to~\citet[Theorem 17.19]{aliprantis1999infinite}, $\Phi_B$ is \lsc{} as a multifunction $(X \times A) \to 2^{B}$.
    \label{rem:a-lsc-implies-lsc}
\end{remark}

In addition, we introduce the following transfer properties of the function $f$. The first property recasts \FPTusc{} of the worst-loss function $f^\#$ in terms of the objective function $f$.

\begin{definition}
    [$B$-uniform feasible path transfer upper semicontinuity]
    Let $E \subset \Dom(\Phi_A)$. The function $f : \Gr(\Phi_B) \subseteq X \times A \times B \to \bar \RR$ is \emph{$B$-uniform feasible path transfer upper semicontinuous} (\emph{$B$-uniform \FPTusc{}}) on $\Gr_{E}(\Phi_B)$ if, for each $(x,a) \in E$ and $\lambda \in \RR$, if $f(x,a,b) + \varepsilon < \lambda$ for each $b \in \Phi_B(x,a)$ and some $\varepsilon > 0$, then there is $\varepsilon' > 0$ and a neighborhood $U(x)$ such that for each $x' \in U(x)$ there exists $a' \in A$ with $(x',a') \in E$ such that $f(x',a',b') + \varepsilon' < \lambda$ for each $b' \in \Phi_B(x',a')$.
    \label{def:b-uniform-fptusc}
\end{definition}

When $X$, $A$, and $B$ are metric spaces, $B$-uniform \FPTusc{} functions admit a sequence characterization, which follows from Theorem~\ref{thm:metric-fptusc}.

\begin{theorem}
    [Characterization of $B$-uniform \FPTusc{} functions for metric spaces $X$, $A$, and $B$]
    Suppose $X$, $A$, and $B$ are metric spaces, and let $x \in X$. The function $f : \Gr(\Phi_B) \subset X \times A \times B \to \bar{\RR}$ is $B$-uniform \FPTusc{} on $\Gr_{\set{x}}(\Phi_B)$ if and only if, for each $a \in \Phi_A(x)$, $\varepsilon > 0$, and sequence $\seq{x_n}_{n=1}^\infty \subset \Dom(\Phi_A)$ with $x_n \to x$, there exist $b \in \Phi_B(x,a)$ and sequences $\seq{a_n}_{n=1}^{\infty}$ and $\seq{b_n}_{n=1}^{\infty}$ such that $a_n \in \Phi_A(x_n)$ and $b_n \in \Phi_B(x_n, a_n)$ for each $n \in \NN$ and such that
    \begin{equation}
        \limsup_{n \to \infty} f(x_n, a_n, b_n) \leq f(x,a,b) + \varepsilon.
        \label{eq:b-uniform-fptusc-metric}
    \end{equation}
    \label{thm:b-uniform-fptusc-metric}
\end{theorem}

\begin{proof}
The inequality
\begin{equation}
    \limsup_{n \to \infty} f^\#(x_n,a_n) \leq f^\#(x,a)
    \label{eq:b-fptusc-metric-ieq1}
\end{equation}
is by definition identical to the inequality
\begin{equation}
    \limsup_{n \to \infty} \sup_{b' \in \Phi_B(x_n, a_n)} f(x_n, a_n,
    b') \leq  \sup_{b \in \Phi_B(x,a)} f(x,a,b).
    \label{eq:b-fptusc-metric-ieq2}
\end{equation}
If, for all $a \in \Phi_A(x)$ and $x_n \to x$, there exists a sequence $a_n \in \Phi(x_n)$ satisfying the inequality in~(\ref{eq:b-fptusc-metric-ieq1}), then it follows that for each $\varepsilon > 0$ there exist $b \in \Phi_B(x,a)$ such that
\begin{equation}
    \limsup_{n \to \infty} \sup_{b' \in \Phi_B(x_n, a_n)} f(x_n, a_n,
    b') \leq f(x,a,b) + \varepsilon.
    \label{eq:b-fptusc-metric-ieq3}
\end{equation}
Moreover, for each $n \in \NN$ one can choose $b_n \in \Phi_B(x_n, a_n)$ such that
\begin{equation*}
    f(x_n, a_n, b_n) + n\inv \geq \sup_{b' \in \Phi_B(x_n, a_n)} f(x_n, a_n, b'),
\end{equation*}
which implies that for each $\varepsilon > 0$ there exists $b \in \Phi_B(x, a)$ and a sequence $\seq{b_n}_{n =1}^\infty$ such that the inequality in~(\ref{eq:b-fptusc-metric-ieq2}) holds. The converse implication follows similarly.
\end{proof}

The following theorem generalizes~\citet[Theorems 5, 6, 7]{feinberg2017continuity}, which gave sufficient conditions for lower semi-continuity, upper semi-continuity, and continuity of the worst-loss function $f^{\#}$. Theorem~\ref{thm:worst-loss-iff-conditions} provides general necessary and sufficient conditions for the worst-loss $f^{\#}$ to be continuous, as well as for existence of solutions in $\Phi^*_B(x,a)$. In particular,

\begin{theorem}
    [Continuity of the worst-loss function]
    Let $f^\# : \Gr(\Phi_A) \subset X \times A \to \bar \RR$ be the worst-loss function defined in Equation (\ref{eq:worst-loss-sup}) for the minimax problem $(X, A, B, f, \Phi_A, \Phi_B)$, and let $x \in X$ and $a \in \Phi_A(x)$. The following statements hold.
    \begin{enumerate}
        \item $f^\#$ is \lsc{} at $(x,a)$ if and only if  $f$ is \FPTlsc{} on $\Gr_{\set{(x,a)}}(\Phi_B)$. \label{item:worst-loss-lsc}
        \item $f^\#$ is \usc{} at $(x,a)$ if and only if $f$ is~\usupsc{} at $(x,a)$. \label{item:worst-loss-usc}
        \item $f^\#$ is \usc{} at $(x,a)$, and the supremum in Equation~(\ref{eq:worst-loss-sup}) can be replaced with a maximum if and only if $f$ is \umaxsc{} at $(x,a)$. \label{item:worst-loss-usc-sol}
        \item $f^\#$ is continuous at $(x,a)$ if and only if $f$ is \FPTlsc{} on $\Gr_{\set{(x,a)}}(\Phi_B)$, and~\usupsc{} at $(x,a)$. \label{item:worst-loss-cont}
        \item $f^\#$ is continuous at $(x,a)$ and the supremum in Equation~(\ref{eq:worst-loss-sup}) can be replaced with a maximum if and only if $f$ is \FPTlsc{} on $\Gr_{(x,a)}(\Phi_B)$ and~\umaxsc{} at $(x,a)$.  \label{item:worst-loss-cont-sol}
    \end{enumerate}
    \label{thm:worst-loss-iff-conditions}
\end{theorem}

\begin{proof}
Consider the minimization problem $(X \times A, B, \Phi_B, -f)$, and observe that the value function is $-f^\#$.  Statement~\ref{item:worst-loss-lsc} follows from Theorem~\ref{thm:general-usc-value}. Statement~\ref{item:worst-loss-usc} follows from Theorem~\ref{thm:linfsc}.  Statement~\ref{item:worst-loss-usc-sol} follows from Theorem~\ref{thm:lminsc}. Statement~\ref{item:worst-loss-cont} follows from combining statements~\ref{item:worst-loss-lsc} and~\ref{item:worst-loss-usc}.  Statement~\ref{item:worst-loss-cont-sol} follows from combining statements~\ref{item:worst-loss-lsc} and~\ref{item:worst-loss-usc-sol}.
\end{proof}

The next theorem provides general sufficient conditions for the solution multifunction $\Phi_B$ to be \usc{} and compact-valued. These conditions also imply that the worst-loss function $f^\#$ is continuous. The theorem generalizes~\citet[Theorem 12]{feinberg2017continuity} by relaxing the assumption that $f$ is continuous.

\begin{theorem}
    [Continuity properties for the solution multifunction $\Phi_B^{\#}$]
    Let $x \in X$ and let $a \in \Phi_A(x)$.  If the function $f$  is \FPTlsc{} on $\Gr_{\set{(x,a)}}(\Phi_A)$ and $f$ is $\KK\NN$-sup-compact on $\Gr_{\set{x}}(\Phi)$, then the worst-loss function $f^\# : \Gr(\Phi_A) \subset X \times A \to \bar \RR$ is continuous at $(x, a)$. Furthermore, if $f^\#(x,a) = +\infty$, then $\Phi_B^{\#}(x,a) = \Phi_B(x,a)$. Otherwise, $\Phi_B^{\#}$ is \usc{} at $(x,a)$, and $\Phi_B^{\#}(x,a)$ is nonempty and compact. \label{thm:complete-worst-loss}
\end{theorem}

\begin{proof}
The results follow from Theorem~\ref{thm:main-berge-max} to the minimization problem $(X \times A, B, \Phi_B, -f)$, where we observe that the value function is $-f^\#$.
\end{proof}

The next properties are equivalent characterizations of \lminsc{} and \linfsc{} of the worst-loss function $f^\#$ in terms of the objective function $f$.

\begin{definition}
    [$B$-feasible path transfer \lminsc{}; $B$-feasible path transfer \linfsc{}]
    Let $E \subset \Dom(\Phi_A)$.
    \begin{enumerate}
        \item The function $f : \Gr(\Phi_B) \subseteq X \times A \times B \to \bar \RR$ is \emph{$B$-feasible path transfer lower min-semicontinuous} (\emph{$B$-\FPTlminsc}) on $\Gr_{E}(\Phi_B)$ if, for each $x \in E$, there exists $a \in \Phi_A(x)$ such that for all $\gamma \in \RR$ with $\gamma < f(x,a,b)$ for some $b \in \Phi_B(x,a)$, there exists a neighborhood $U(x)$ such that for each $x' \in U(x)$ and $a' \in \Phi_A(x')$, there exists $b' \in \Phi_B(x',a')$ such that $\gamma < f(x', a', b')$.
        \item The function $f$ is \emph{$B$-feasible path transfer lower inf-semicontinuous} (\emph{$B$-\FPTlinfsc{}}) on $\Gr_E(\Phi_B)$ if, for each $x \in E$ and $\gamma \in \RR$, if for each $a \in \Phi_A(x)$ there exists $b \in \Phi_B(x,a)$ such that $\gamma < f(x,a,b)$, then there exists a neighborhood $U(x)$ such that for each $x' \in U(x)$ and $a' \in \Phi_A(x')$, there exists $b' \in \Phi_B(x',a')$ such that $\gamma < f(x',a',b')$.
    \end{enumerate}
    \label{def:bfptlinfsc-bfptlinfsc}
\end{definition}

The following theorem is a generalization of~\citet[Theorems 8, 9, 10]{feinberg2017continuity}, which gives sufficient conditions for lower semicontinuity, upper semicontinuity, and continuity of the minimax function $f^*$, respectively.  Theorem~\ref{thm:minimax-iff-conditions} provides general necessary and sufficient conditions for the minimax function $f^*$ to be continuous and for $\Phi^*_A(x) \ne \emptyset$.

\begin{theorem}
    [Continuity of the minimax function]
    For $x \in X$ and for the minimax problem $(X, A, B, f, \Phi_A, \Phi_B)$, the following statements hold.
    \begin{enumerate}
        \item The minimax function $f^* : X \to \bar \RR$ is \usc{} at $x$ if and only if the objective function $f : \Gr(\Phi_B) \subset X \times A \times B \to \bar \RR$ is $B$-uniform \FPTusc{} on $\Gr_{\Gr_{\set{x}}}(\Phi_B)$; \label{item:minimax-usc}
        \item $f^*$ is \lsc{} at $x$ if and only if $f$ is $B$-\FPTlinfsc{} on $\Gr_{\Gr_{\set{x}}(\Phi_A)}(\Phi_B)$; \label{item:minimax-lsc}
        \item $f^*$ is \lsc{} at $x$, and the infimum in Equation~(\ref{eq:minimax-inf}) can be replaced with a minimum if and only if $f$ is $B$-\FPTlminsc{} on $\Gr_{\Gr_{\set{x}}(\Phi_A)}(\Phi_B)$; \label{item:minimax-lsc-sol}
        \item $f^*$ is continuous at $x$ if and only if $f$ is $B$-uniform \FPTusc{} and $B$-\FPTlinfsc{} on $\Gr_{\set{x}(\Phi_A)}(\Phi_B)$; \label{item:minimax-cont}
        \item $f^*$ is continuous at $x$, and the infimum in Equation~(\ref{eq:minimax-inf}) can be replaced with a minimum if and only if $f$ is $B$-uniform \FPTusc{}, and $B$-\FPTlminsc{} on $\Gr_{\set{x}(\Phi_A)}(\Phi_B)$.  \label{item:minimax-cont-sol}
    \end{enumerate}
    \label{thm:minimax-iff-conditions}
\end{theorem}

\begin{proof}
Statement~\ref{item:minimax-cont} follows by combining statements~\ref{item:minimax-usc} and~\ref{item:minimax-lsc}. Statement~\ref{item:minimax-cont-sol} follows by combining statements~\ref{item:minimax-usc} and~\ref{item:minimax-lsc-sol}.

To prove statement~\ref{item:minimax-usc}, in view of Theorem~\ref{thm:general-usc-value}, it suffices to show that $f$ is $B$-uniform \FPTusc{} on $\Gr_{\set{x}(\Phi_A)}(\Phi_B)$ if and only if the worst-loss function $f^{\#} : \Gr(\Phi_A) \subset X \times A \to \bar \RR$ is \FPTusc{} on $\Gr_{\set{x}}(\Phi_A)$.  Assume first that $f$ is $B$-uniform \FPTusc{}.  Let $a \in \Phi_A(x)$, and suppose $\lambda \in \RR$ with $f^{\#}(x,a) < \lambda$. Let $\lambda' \in (f^{\#}(x,a), \lambda)$ be arbitrary, and take $\varepsilon = \lambda - \lambda'$. Then for each $b \in \Phi_B(x,a)$, there holds $f(x,a,b) + \varepsilon < \lambda$, so we can fix $\varepsilon' > 0$ and a neighborhood $U(x)$ such that for each $x' \in U(x)$ there exists $a' \in \Phi_A(x')$ such that $f(x',a',b') + \varepsilon' < \lambda$ for each $b' \in \Phi_B(x',a')$.  Then $f^{\#}(x',a') + \varepsilon' = \sup_{b' \in \Phi_B(x',a')} f(x',a',b') + \varepsilon' \leq  \lambda$, which implies $f^{\#}(x',a') < \lambda$, as needed.  Conversely, let $a \in \Phi_A(x)$, and suppose $\lambda \in \RR$ satisfies $f(x,a,b) + \varepsilon < \lambda$ for some $\varepsilon > 0$ and each $b \in \Phi_B(x,a)$. Then $f^{\#}(x,a) < \lambda$, so there is a neighborhood $U(x)$ such that for each $x' \in U(x)$ there is $a' \in \Phi_A(x')$ such that $f^{\#}(x',a') < \lambda$. Taking $\varepsilon' = \min(1, \lambda - f^{\#}(x',a'))$ implies that $f(x',a',b') + \varepsilon' < \lambda$ for each $b' \in \Phi_B(x',a')$, as needed.

Next, to prove statement~\ref{item:minimax-lsc}, in view of Theorem~\ref{thm:linfsc}, it suffices to show that $f$ is $B$-\FPTlinfsc{} on $\Gr_{\set{x}(\Phi_A)}(\Phi_B)$ if and only if $f^{\#}$ is \linfsc{} at $x$. Assume that $f$ is $B$-\FPTlinfsc{}. Let $\gamma \in \RR$, and suppose for all $a \in \Phi_A(x)$ one has $\gamma < f^{\#}(x,a)$.  Then for each $a \in \Phi_A(x)$ there exists $b \in \Phi_B(x,a)$ such that $\gamma < f(x,a,b)$. Therefore, there is a neighborhood $U(x)$ such that for each $x' \in U(x)$ and $a' \in \Phi_A(x')$ there exists $b' \in \Phi_B(x',a')$ such that $\gamma < f(x', a', b') \leq f^{\#}(x',a')$, establishing that $f^{\#}$ is \linfsc{}. Conversely, if $f^{\#}$ is \linfsc{}, let $\gamma \in \RR$ and suppose that for each $a \in \Phi_A(x)$ there exists $b \in \Phi_B(x,a)$ such that $\gamma < f(x,a,b)$. Then $\gamma < f^{\#}(x,a)$, so there is a neighborhood $U(x)$ such that for each $x' \in U(x)$ and $a' \in \Phi_A(x')$, one has $\gamma < f^{\#}(x',a')$. This implies there exists $b' \in \Phi_B(x',a')$ such that $\gamma < f(x',a',b')$, establishing that $f$ is $B$-\FPTlinfsc{}.

Finally, to prove statement~\ref{item:minimax-lsc-sol}, in view of Theorem~\ref{thm:lminsc}, it suffices to show that $f$ is $B$-\FPTlminsc{} on $\Gr_{\set{x}(\Phi_A)}(\Phi_B)$ if and only if $f^{\#}$ is \lminsc{} at $x$. Assume first that $f$ is $B$-\FPTlminsc{}. Fix $a^* \in \Phi_A(x)$ as specified in Definition~\ref{def:bfptlinfsc-bfptlinfsc}.  Then for all $\gamma \in \RR$ such that $\gamma < f^{\#}(x,a^*)$, let $\gamma' \in (\gamma, f^{\#}(x,a^*))$.  Then by definition of $f^{\#}$, there exists $b \in \Phi_B(x,a^*)$ such that $\gamma' < f(x,a^*,b)$. Then we can fix a neighborhood $U(x)$ such that for each $x' \in U(x)$ and $a' \in \Phi_A(x')$, there exists $b' \in \Phi_B(x',a')$ such that $\gamma' < f(x',a',b')$; hence, $\gamma < \gamma' \leq f^{\#}(x',a')$, establishing that $f^{\#}$ is \lminsc{}.  Conversely, if $f^{\#}$ is \lminsc{}, let $x \in X$ and $\gamma \in \RR$. Fix $a^* \in \Phi_A(x)$ as in Definition~\ref{def:lminsc}. Then if there is $b \in \Phi_B(x,a^*)$ such that $\gamma < f(x,a^*,b)$, observe that $\gamma < f(x,a^*,b) \leq f^{\#}(x,a^*)$, so there exists a neighborhood $U(x)$ such that for each $x' \in U(x)$ and $a' \in \Phi_A(x')$ one has $\gamma < f^{\#}(x',a')$; hence, there exists $b' \in \Phi_B(x',a')$ such that $\gamma < f(x',a',b')$, establishing that $f$ is $B$-\FPTlminsc{}.
\end{proof}

The next theorem states general sufficient conditions for the solution multifunction $\Phi_A^*$ to be \usc{} and compact-valued. This generalizes Theorem 13 in~\citet{feinberg2017continuity} by relaxing the assumption that $f$ is continuous.

\begin{theorem}
    [Continuity properties of the solution multifunction $\Phi_A^*$]
    Suppose the multifunction $\Phi_B : \Gr(\Phi_A) \subseteq X \times A \to 2^B$ is $A$-lower semicontinuous on $\Gr_{\Gr_{\set{x}}(\Phi_A)}(\Phi_B)$, the function $f^{A \leftrightarrow B} : \Gr(\Phi_B^{A \leftrightarrow B}) \subseteq X \times B \times A \to \bar \RR$ is $\KK\NN$-inf-compact on $\Gr_{\Gr_{\set{x}}(\Phi_A^{A \leftrightarrow B})}(\Phi_B^{A \leftrightarrow B})$, and the function $f : \Gr(\Phi_B) \subseteq X \times A \times B \to \bar \RR$ is $B$-uniform \FPTusc{} on $\Gr_{\Gr_{\set{x}}(\Phi_A)}(\Phi_B)$.  Then the minimax function $f^* : X \to \bar \RR$ is continuous at $x$, and moreover,
    \begin{enumerate}
        \item the infimum in Equation~(\ref{eq:minimax-inf}) can be replaced with a minimum;
        \item if $f^*(x) < +\infty$, then $\Phi_A^*$ is upper semicontinuous at $x$, and $\Phi_A^*(x)$ is nonempty and compact.
    \end{enumerate}
    \label{thm:complete-minimax}
\end{theorem}

\begin{proof}
    The function $f^*$ is \usc{} in view of Theorem~\ref{thm:minimax-iff-conditions}\ref{item:minimax-usc}. From Remark~\ref{rem:a-lsc-implies-lsc} in view of Theorem~\ref{thm:sufficient-usc}, the function $-f$ is \FPTusc{}, and by Theorem~\ref{thm:general-usc-value}, $-f^{\#}$ is \usc{}, so $f^{\#}$ is \lsc{}. Next, suppose $\seq{(x_\alpha, a_\alpha)}_{\alpha \in I} \subseteq \Gr(\Phi_A)$ with $x_\alpha \to x$ and $\sup_{\alpha \in I} f^{\#}(x_\alpha, a_\alpha) < +\infty$. Let $\bar{a} \in \Phi_A(x)$ and  $b \in \Phi_B(x,\bar{a})$ be arbitrary. Since $\Phi_B$ is $A$-lower semicontinuous, there exists a net $\seq{b_\alpha}_{\alpha \in I}$ such that $b_\alpha \in \Phi_B(x_\alpha, a_\alpha)$ for each $\alpha$ and a subnet $\seq{b_\beta}_{\beta \in J} \subseteq \seq{b_\alpha}_{\alpha \in I}$ such that $b_\beta \to b$. Since
    \begin{equation*}
        \sup_{\beta \in J} f^{A \leftrightarrow B}(x_\beta, b_\beta, a_\beta)
        \leq \sup_{\alpha \in I} f^{\#}(x_\alpha, a_\alpha) < +\infty,
    \end{equation*}
    and since $f^{A \leftrightarrow B}$ is $\KK\NN$-inf-compact, there exists a subnet $\seq{a_\gamma}_{\gamma \in L} \subseteq \seq{a_\beta}_{\beta \in J}$ such that $a_\gamma \to a \in \Phi_B^{A \leftrightarrow B}(x,b)$. Since $a \in \Phi_A(x)$, this proves that $f^{\#}$ is $\KK\NN$-inf-compact on $\Gr_{\set{x}}(\Phi_A)$. Since $f^{\#}$ is \FPTusc{} and $\KK\NN$-inf-compact on $\Gr_{\set{x}}(\Phi_A)$, the results follow from applying Theorem~\ref{thm:main-berge-max}.
\end{proof}

Both of the above theorems are local formulations of the minimax problem. This is because continuity of $f^*$ and of $\Phi_A^*$ are themselves local properties. Of course, the theorems can each be naturally reformulated as global versions. For example, we have the following as a corollary of Theorem~\ref{thm:complete-minimax}.
\begin{corollary}
    Suppose the multifunction $\Phi_B : \Gr(\Phi_A) \subseteq X \times A \to 2^B$ is $A$-lower semicontinuous on $\Gr(\Phi_B)$, the function $f^{A \leftrightarrow B} : \Gr(\Phi_B^{A \leftrightarrow B}) \subseteq X \times B \times A \to \bar \RR$ is $\KK\NN$-inf-compact on $\Gr(\Phi_B^{A \leftrightarrow B})$, and the function $f : \Gr(\Phi_B) \subseteq X \times A \times B \to \bar \RR$ is $B$-uniform \FPTusc{} on $\Gr(\Phi_B)$.  Then the minimax function $f^* : X \to \bar \RR$ is continuous on $\Dom(\Phi_A)$, and for each $x \in \Dom(\Phi_A)$,
    \begin{enumerate}
        \item the infimum in Equation~(\ref{eq:minimax-inf}) can be replaced with the minimum;
        \item if $f^*(x) < +\infty$, then $\Phi_A^*$ is upper semicontinuous at $x$ and $\Phi_A^*(x)$ is nonempty and compact.
    \end{enumerate}
    \label{thm:global-complete-minimax}
\end{corollary}

\begin{proof}
For each $x \in \Dom(\Phi_A)$, the multifunction $\Phi_B$ is $A$-lower semicontinuous on $\Gr_{\Gr_{\set{x}}(\Phi_A)}(\Phi_B)$, the function $f^{A \leftrightarrow B}$ is $\KK\NN$-inf-compact on $\Gr_{\Gr_{\set{x}}(\Phi_A^{A \leftrightarrow B})}(\Phi_B^{A \leftrightarrow B})$, and the function $f$ is $B$-uniform \FPTusc{} on $\Gr_{\Gr_{\set{x}}(\Phi_A)}(\Phi_B)$. The result thus follows from Theorem~\ref{thm:complete-minimax}.
\end{proof}

In Section~\ref{sec:inventory-control} we considered applications of Theorem~\ref{thm:main-berge-max} to inventory control problems with setup costs, possible lost sales, possibly unbounded orders, and possibly unbounded storage capacity. A natural direction following the results of this section, in particular Theorem~\ref{thm:complete-minimax}, is to consider robust optimization of the aforementioned problems. Although such results are beyond the scope of the present work, they can be applied to robust inventory control, which is a rich field of study; see e.g.,~\citet{bertsimas2006robust} and~\citet{thorsen2017robust}.
\section*{Acknowledgements}

The first author thanks R. Tyrrell Rockafellar and Johannes O. Royset for valuable comments. The second author was partially supported by the National Research Foundation of Ukraine, Grant No. 2020.01/0283.

\bibliography{FeinbergKasyanovKraemer}

\end{document}